\documentclass[12pt]{amsart}
\usepackage{ben}

\title[Non-autonomous parabolic Cauchy problems]{On well-posedness for non-autonomous parabolic Cauchy problems with rough initial data}
\author{Hedong Hou}
\address{Universit{\'e} Paris-Saclay, CNRS, Laboratoire de Math\'{e}matiques d'Orsay, 91405 Orsay, France}
\email{hedong.hou@universite-paris-saclay.fr}

\date{May 14, 2025}

\begin{document}

\begin{abstract}
    We establish a complete picture for existence, uniqueness, and representation of weak solutions to non-autonomous parabolic Cauchy problems of divergence type. The coefficients are only assumed to be uniformly elliptic, bounded, measurable, and complex-valued, without any additional regularity or symmetry conditions. The initial data are tempered distributions taken in homogeneous Hardy--Sobolev spaces $\dot{H}^{s,p}$, and source terms belong to certain scales of weighted tent spaces. Weak solutions are constructed with their gradients in weighted tent spaces $T^{p}_{s/2}$. Analogous results are also exhibited for initial data in homogeneous Besov spaces $\dot{B}^{s}_{p,p}$.
\end{abstract}
\maketitle
\tableofcontents
\section{Introduction}
\label{sec:intro}

The objective of this paper is to present a complete picture for existence, uniqueness, and representation of weak solutions to the non-autonomous parabolic Cauchy problem of divergence type
\begin{equation}
    \label{e:NaPC}
    \begin{cases}
        \partial_t u - \Div_x( A(t,x) \nabla_x u ) = f + \Div_x F, \quad (t,x) \in (0,\infty) \times \bR^n =: \bR^{1+n}_+ \\
        u(0)=u_0
    \end{cases}.
\end{equation}
By the initial condition $u(0)=u_0$, we require $u(t)$ tends to 0 as $t \to 0$ in distributional sense, \textit{i.e.}, in $\scrD'(\bR^n)$.

Assume that the coefficient matrix $A \in L^\infty(\bR^{1+n}_+;\mat_n(\bC))$ is \emph{uniformly elliptic}. Namely, there exist $\Lambda_0, \Lambda_1>0$ so that for a.e. $(t,x) \in (0,\infty) \times \bR^n$ and any $\xi,\eta \in \bC^n$,
\begin{equation}
    \label{e:ellpitic}
    \Re(\langle A(t,x)\xi,\xi \rangle) \geq \Lambda_0 |\xi|^2, \quad |\langle A(t,x)\xi,\eta \rangle| \leq \Lambda_1 |\xi| |\eta|.
\end{equation}
It is worth pointing out that we do \emph{not} impose \emph{any} assumptions on symmetry or regularity of the coefficients in either the time variable $t$ or the space variable $x$. 

This model serves as a simplified but representative example for parabolic systems. In fact, our results presented herein readily extend to parabolic systems where the ellipticity condition in \eqref{e:ellpitic} is replaced by a Gårding inequality on $\bR^n$ uniformly with respect to $t$. That is, there exists $\Lambda_0>0$ so that for a.e. $t>0$ and any $\varphi \in \DotH^1(\bR^n;\bC^m)$,
\begin{equation}
    \label{e:Garding-elliptic}
    \Re \int_{\bR^n} (A(t,x) \nabla \varphi(x)) \cdot \overline{\nabla \varphi}(x) \ge \Lambda_0 \|\nabla \varphi\|_{L^2(\bR^n)}^2.
\end{equation}
Here, $A(t,x) = (A^{\alpha,\beta}_{i,j}(t,x))^{1 \le \alpha,\beta \le n}_{1 \le i,j \le m}$ belongs to $L^\infty(\bR^{1+n}_+;\mat_{mn}(\bC))$, and we use the notation
\[ (A(t,x) \nabla \varphi(x)) \cdot \overline{\nabla \varphi}(x) := \sum_{\substack{1 \le i,j \le m \\ 1 \le \alpha,\beta \le n}} A^{\alpha,\beta}_{i,j}(t,x) \partial_\beta \varphi^j(x) \overline{\partial_\alpha \varphi^i}(x).  \]
The homogeneous Sobolev space $\DotH^1(\bR^n;\bC^m)$ is the closure of $\Cc(\bR^n;\bC^m)$ with respect to the semi-norm $\|\nabla \varphi\|_{L^2(\bR^n;\bC^m)}$.

We show that for nearly optimal ranges of parameters $s$ and $p$, given initial data $u_0$ as a tempered distribution in the \emph{homogeneous Hardy--Sobolev space} $\DotH^{s,p}$ and suitable source terms $f$ and $F$, one can construct weak solutions $u$ to \eqref{e:NaPC}. The solution class is characterized by the property that $\nabla u$ belongs to the weighted tent space $T^p_{s/2}$. Uniqueness and representation of weak solutions in this solution class are also established. Precise definitions of theses function spaces are given in Section \ref{sec:spaces}. For context, the spaces $\DotH^{s,p}$ are defined via Littlewood--Paley decomposition, where $s$ denotes regularity and $p$ denotes integrability. When $s=0$, $\DotH^{0,p}$ identifies with the Hardy space $H^p$ if $0<p \le 1$, the Lebesgue space $L^p$ if $1<p<\infty$, and the John--Nirenberg bounded mean oscillation space $\bmo$ if $p=\infty$.

A classical approach to this problem involves establishing \textit{a priori} mixed $L^p_t L^q_x$-estimates for weak solutions and their derivatives, particularly the highest order derivative $\nabla u$. However, to the best of our knowledge, obtaining such estimates for arbitrary coefficient matrices $A$ remains an open problem. Partial results exist under additional regularity assumptions. For instance, when $x \mapsto A(t,x)$ is uniformly continuous, the problem has been addressed by Lady{\v z}enskaja, Solonnikov and Ural'ceva \cite{Ladyzhenskaya-Solonnikov-Uralceva1968-parabolic-PDE}. Krylov \cite{Krylov2007-Lp-div-VMOx} later relaxes this condition, requiring only that $x \mapsto A(t,x)$ belongs to $\vmo$, see also the works of Dong and Kim on parabolic systems \cite{Dong-Kim2011-Lp-div-VMOx-system,Dong-Kim2011-Lp-VMOx-system,Dong-Kim2018-Lp-pe-Ap-weight}. 

The fundamental strategy underlying these results is to first derive the estimates for simplified cases, such as constant or space-independent coefficients, before extending them to more general ones via perturbation arguments, for which the regularity assumptions on coefficients seem to be unavoidable.

Instead, we consider an alternative framework by using \emph{(weighted) tent spaces} $T^p_\beta$, which leads to finer estimates for rough coefficients and rough initial data. Tent spaces were introduced by Coifman, Meyer and Stein \cite{Coifman-Meyer-Stein1985_TentSpaces}, and originated from the work of Fefferman and Stein \cite{Fefferman-Stein1972Hp} on real Hardy spaces. The key insight is that their norms first capture localized $L^2$ integrability in the interior (in both time and space) and then incorporate an additional control (with weight in time) to access the limit at the boundary (materialized as $\bR^n$), especially where standard mixed $L^p_t L^q_x$-spaces fall short.

In order to understand the role of tent spaces, we first consider the homogeneous Cauchy problem
\begin{equation*}
    \tag{HC}
    \begin{cases}
        \partial_t u - \Div_x(A(t,x)\nabla_x u)=0, \quad 0<t<T, ~ x \in \bR^n \\
        u(0)=u_0
    \end{cases},
\end{equation*}
where the complex-valued coefficient matrix $A$ is only assumed to be uniformly elliptic, bounded, and measurable. The fundamental work of Lions \cite{Lions1957_L2} first establishes well-posedness of weak solutions to \eqref{e:hc} with initial data $u_0 \in L^2$ in the energy space $L^\infty((0,T);L^2) \cap L^2((0,T);W^{1,2})$ for finite time $T>0$. In the special case of real coefficients, Aronson \cite{Aronson1968-Propagator} constructs fundamental solutions $K(t,s;x,y)$ with pointwise Gaussian decay and shows that the \emph{propagators} $(\Gamma_A(t,s))_{0 \le s < t < T}$ defined on $L^2(\bR^n)$ by
\[ (\Gamma_A(t,s) u_0)(x) := \int_{\bR^n} K(t,s;x,y) u_0 (y) dy \]
represent the unique weak solution $u(t,x)$ in the energy space to
\begin{equation}
    \label{e:hc-(s,T)}
    \begin{cases}
        \partial_t u - \Div_x(A(t,x)\nabla_x u)=0, \quad s<t<T, ~ x \in \bR^n \\
        u(s)=u_0
    \end{cases}.
\end{equation}

A recent work of Auscher, Monniaux and Portal \cite{Auscher-Monniaux-Portal2019Lp} imports the tent space theory to establish well-posedness of weak solutions to \eqref{e:hc} with initial data $u_0 \in L^p(\bR^n)$ for $2-\epsilon<p<\infty$, where $\epsilon>0$ is a constant depending on the ellipticity of $A$. To this end, they first extend Lions' $L^2$-theory to $T=\infty$ within the solution class $\nabla u \in L^2(\bR^{1+n}_+)$. Note that $L^2(\bR^{1+n}_+)$ identifies with the tent space $T^2_0$. By shifting the time, they define the propagators $(\Gamma_A(t,s))_{0 \le s \le t < \infty}$ as a family of contractions on $L^2(\bR^n)$ so that for any $s \ge 0$ and $u_0 \in L^2$,
\begin{equation}
    \label{e:def-propagator}
    u(t,x) := (\Gamma_A(t,s) u_0)(x)
\end{equation}
is the unique weak solution to \eqref{e:hc-(s,T)} with $\nabla u \in L^2((s,\infty) \times \bR^n)$. This definition agrees with Aronson's definition by fundamental solutions for real coefficients but there might not be a representation with a kernel.

The weak solution $u$ is constructed by extension of the \emph{propagator solution map} $\cE_A$, initially defined from $L^2$ to $L^2_{\loc}((0,\infty);W^{1,2}_{\loc})$ by
\begin{equation}
    \label{e:EA}
    \cE_A(u_0)(t,x) := (\Gamma_A(t,0)u_0)(x), \quad t>0, ~ x \in \bR^n.
\end{equation}
The solution class is characterized by $\nabla u \in T^p_0$, with the equivalence
\[ \|\nabla \cE_A(u_0)\|_{T^p_0} \eqsim \|u_0\|_{L^p}. \]
A later complement by Zato\'n \cite{Zaton2020wp} extends it to $p=\infty$ with
\[ \|\nabla \cE_A(u_0)\|_{T^\infty_0} \eqsim \|u_0\|_{\bmo}. \]
Zato\'n also proves uniqueness and representation of weak solutions in this solution class and extends these results to parabolic systems.

Motivated by applications towards stochastic analysis, Portal and Veraar propose the following problem on well-posedness of the homogeneous Cauchy problem \eqref{e:hc} with initial data in $\DotH^{s,p}$.
\begin{prob}[{\cite[p.583]{Portal-Veraar2019_Stochastic-MR-tent}}]
    \label{prob:PV-eq}
    Let $A \in L^\infty(\bR^{1+n}_+; \mat_n(\bC))$ be uniformly elliptic. Does there exist a range of $s \ge 0$ and $p \in (1,\infty)$ so that the equivalence
    \begin{equation}
        \label{e:PV-u0-nabla-u-eq}
         \|\nabla \cE_A(u_0)\|_{T^p_{s/2}} \eqsim \|u_0\|_{\DotH^{s,p}}
    \end{equation}
    holds for any $u_0 \in \DotH^{s,p}$?
\end{prob}

The equivalence \eqref{e:PV-u0-nabla-u-eq} clearly relates regularity of the initial data to the time weight in the tent space norm of the solution. In the special case $A=\bI$ (\textit{i.e.}, $\Div(A\nabla)$ agrees with the Laplacian $\Delta$), they intend to prove it for all $s \ge 0$ and $p \in (1,\infty)$, but the proof has a gap. In fact, we have shown in \cite[Theorem 1.1]{Auscher-Hou2024-HCL} that for $s \ge 1$ and $1 \le p \le \infty$, the equivalence \eqref{e:PV-u0-nabla-u-eq} holds if and only if $u_0$ is constant, for which $\|u_0\|_{\DotH^{s,p}}=0$.

More generally, our previous work \cite{Auscher-Hou2024-HCL} addresses the problem for all time-independent coefficients $A=A(x)$. Indeed, the theory of maximal accretive operators yields that $-L:=\Div(A\nabla)$ generates an analytic semigroup $(e^{-tL})_{t \ge 0}$ on $L^2$, so the propagators become
\[ \boxed{\Gamma_A(t,s) = e^{-(t-s)L}, \quad \text{ if } A=A(x)}. \]
We show that for $-1<s<1$ and $\tilp_L(s)< p \le \infty$, the solution map $u_0 \mapsto u(t,x)=(e^{-tL} u_0)(x)$ is an \emph{isomorphism} from $\DotH^{s,p}+\bC$ to the space of weak solutions to $\partial_t u - \Div (A(x)\nabla u) = 0$ with $\nabla u \in T^p_{s/2}$. In particular, this yields the equivalence \eqref{e:PV-u0-nabla-u-eq}. The lower bound $\tilp_L(s)$ has an explicit expression depending on $s$ and $L$. For the identity matrix $A=\bI$ ($L=-\Delta$), we have $\tilp_{-\Delta}(s) = \frac{n}{n+s+1}<1$.

In this paper, we address Problem \ref{prob:PV-eq}. We also consider source terms in tent spaces. To get \textit{a priori} estimates for the tent space norms of $u$ and $\nabla u$, we employ the theory of singular integral operators on tent spaces developed in \cite{Auscher-Hou2025_WPMRTent}. A crucial component is the $L^p-L^q$ off-diagonal decay of the propagators (see Definition \ref{def:odd}). To prove uniqueness, we use an internal representation of weak solutions introduced in \cite{Auscher-Monniaux-Portal2019Lp}. Notably, for $s=0$, our approach not only recovers Zato\'n's uniqueness results via a conceptually simpler proof, but also relaxes the pointwise elliptic condition, see Remark \ref{rem:rep-system-ellp}.

\subsection{Main results}
\label{ssec:intro-results}
For simplicity, we state our results for a single complex-valued equation and leave the detailed formulation for parabolic systems to the reader.

Let $A \in L^\infty(\bR^{1+n}_+;\mat_n(\bC))$ be uniformly elliptic (see \eqref{e:ellpitic}) and $(\Gamma_A(t,s))_{0 \le s \le t < \infty}$ be the propagators associated to $A$ (see \eqref{e:def-propagator}). Let $p_\pm(A)$ be two real numbers in $[1,\infty]$ so that $(p_-(A),p_+(A))$ is the largest open set of exponents $p$ for which $(\Gamma_A(t,s))$ extends to a uniformly bounded family on $L^p$. It is known that $p_-(A)<2$ and $p_+(A)>2$, see Proposition \ref{prop:Gamma-Lp}. For time-independent coefficients, they agree with the critical numbers $p_\pm(L)$ introduced by \cite[\S 3.2]{Auscher2007memoire}. In particular, for the identity matrix $A=\bI$, we have $p_-(\bI) = 1$ and $p_+(\bI) = \infty$.

To precise our results, we introduce several critical exponents. For convenience, we parametrize our results by $\beta \in \bR$, whose relation with the regularity index $s$ is given by
\[ \boxed{s=2\beta+1}. \]

\begin{figure}[ht]
    \centering
    \includegraphics[width=0.6\linewidth]{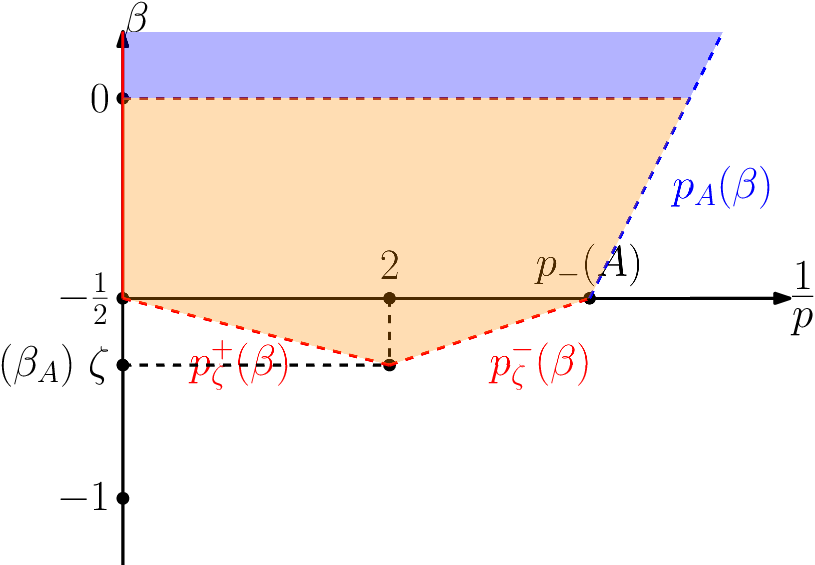}
    \caption{Region of well-posedness}
    \label{fig:exponents}
\end{figure}

Define the number
\begin{equation}
    \label{e:pA(beta)}
    p_A(\beta) := \frac{np_-(A)}{n+(2\beta+1)p_-(A)}.
\end{equation}
Let $\zeta<-1/2$ be a fixed reference number. Define the critical exponents $p^\pm_\zeta(\beta)$ by
\begin{equation}
    \label{e:p-zeta(beta)}
    p^-_\zeta(\beta) := 
    \begin{cases}
        \frac{2(2\zeta+1)p_-(A)}{4(\zeta-\beta)+(2\beta+1)p_-(A)} & \text{ if } \zeta \le \beta<-1/2 \\
        p_A(\beta) & \text{ if } \beta \ge -1/2
    \end{cases},
\end{equation}
and
\begin{equation}
    \label{e:p+zeta(beta)}
    p^+_\zeta(\beta) := 
    \begin{cases}
        \frac{2(2\zeta+1)}{2\beta+1} & \text{ if } \zeta \le \beta<-1/2 \\
        +\infty & \text{ if } \beta \ge -1/2
    \end{cases}.
\end{equation}
Note that $p^-_\zeta(\zeta) = p^+_\zeta(\zeta)=2$. To illustrate these exponents, we give graphic representations in Figure \ref{fig:exponents}. In this figure, we write $p$ for $1/p$ to ease the presentation. We use red dashed lines for graphs of $p^{\pm}_\zeta(\beta)$ when $\beta<-1/2$, red normal line for that of $p^+_\zeta(\beta)$ when $\beta \ge -1/2$, and blue dashed line for that of $p_A(\beta)$. Parallel lines to the blue dashed line are lines of embedding for homogeneous Hardy--Sobolev spaces and weighted tent spaces going downward.

We shall introduce a new parameter $\beta_A \in [-1,-1/2)$, which only depends on the ellipticity of $A$ and the dimension $n$, as the lower bound of $\beta$. Taking $\zeta=\beta_A$, the orange shaded trapezoid becomes the region of well-posedness for initial data in $\DotH^{2\beta+1,p}$, while the blue shaded is for constant initial data. Note that the extreme point $(0,-1/2)$ is always included.

Our first theorem is devoted to the homogeneous Cauchy problem ($f=0$, $F=0$).
\begin{theorem}[Well-posedness of the homogeneous Cauchy problem]
    \label{thm:wp-hc}
    There exists $\beta_A \in [-1,-1/2)$ only depending on the ellipticity of $A$ and the dimension $n$ so that when $\beta_A<\beta<0$ and
    \[ \begin{cases}
        p^-_{\beta_A}(\beta) < p \le p^+_{\beta_A}(\beta) = \infty & \text{ if } \beta \ge -1/2 \\
        p^-_{\beta_A}(\beta) < p < p^+_{\beta_A}(\beta) & \text{ if } \beta_A < \beta < -1/2 
    \end{cases}, \]
    for any $u_0 \in \DotH^{2\beta+1,p}$, there is a unique global weak solution $u$ to the \emph{homogeneous Cauchy problem}
    \begin{equation*}
        \tag{HC}
        \begin{cases}
            \partial_t u - \Div_x(A(t,x)\nabla_x u)=0, \quad (t,x) \in (0,\infty) \times \bR^n \\
            u(0)=u_0
        \end{cases}
    \end{equation*}
    with $\nabla u \in T^p_{\beta+1/2}$. The equivalence holds that
    \[ \|\nabla u\|_{T^p_{\beta+1/2}} \eqsim \|u_0\|_{\DotH^{2\beta+1,p}}. \]
    Moreover, $u$ belongs to $C([0,\infty);\scrS')$.
\end{theorem}

Weak solutions are constructed by extension of the propagator solution map $\cE_A$ (see \eqref{e:EA}) to $\DotH^{s,p}$. See Theorem \ref{thm:hc} for more properties.

The constraint $\beta<0$ ($s<1$) is sharp. In fact, we establish a converse statement, which asserts that any global weak solution $u$ to the equation $\partial_t u - \Div(A\nabla u) = 0$ with $\nabla u \in T^p_{\beta+1/2}$ is uniquely determined by its distributional limit of $u(t)$ as $t \to 0$, called \emph{trace}. When $\beta \ge 0$ ($s \ge 1$), the trace is constant, and so is the solution itself.
\begin{theorem}[Representation]
    \label{thm:rep}
    Let $\beta>\beta_A$ and
    \[ \begin{cases}
        p^-_{\beta_A}(\beta) < p \le p^+_{\beta_A}(\beta) = \infty & \text{ if } \beta \ge -1/2 \\
        p^-_{\beta_A}(\beta) < p < p^+_{\beta_A}(\beta) & \text{ if } \beta_A < \beta < -1/2 
    \end{cases}. \]
    Let $u$ be a global weak solution to $\partial_t u - \Div (A\nabla u) = 0$ with $\nabla u \in T^p_{\beta+{1}/{2}}$. Then $u$ has a trace $u_0 \in \scrS'$, in the sense that $u(t)$ converges to $u_0$ in $\scrS'$ as $t \to 0$. Moreover,
    \begin{enumerate}[label=\normalfont(\roman*)]
        \item 
        \label{item:rep_beta>0}
        If $\beta \ge 0$ and $\frac{n}{n+2\beta} \le p \le \infty$, then $u$ is a constant.
        
        \item 
        \label{item:rep_-1<beta<0}
        If $\beta_A<\beta<0$, then there exist $g \in \DotH^{2\beta+1,p}$ and $c \in \bC$ so that $u_0=g+c$ and $u=\cE_A(g)+c$, where $\cE_A$ is the extension of the propagator solution map defined by \eqref{e:EA}.
    \end{enumerate}
\end{theorem}

Combining Theorems \ref{thm:wp-hc} and \ref{thm:rep} gives a complete answer to Problem \ref{prob:PV-eq}, even with extra ranges for $s<0$ and $p \le 1$. Moreover, it yields the following stronger result. For convenience, define $\cE_A(c)=c$ for any constant function $c$.
\begin{cor}[Isomorphism]
    \label{cor:iso}
    Let $\beta_A<\beta<0$ and 
    \[ \begin{cases}
        p^-_{\beta_A}(\beta) < p \le p^+_{\beta_A}(\beta) = \infty & \text{ if } \beta \ge -1/2 \\
        p^-_{\beta_A}(\beta) < p < p^+_{\beta_A}(\beta) & \text{ if } \beta_A < \beta < -1/2 
    \end{cases}. \]
    The propagator solution map $u_0 \mapsto u=\cE_A(u_0)$ is an isomorphism from $\DotH^{2\beta+1,p} + \bC$ to the space of global weak solutions $u$ to $\partial_t u - \Div (A\nabla u) = 0$ with $\nabla u \in  T^p_{\beta+1/{2}}$, and 
    \[ \|\nabla u\|_{T^p_{\beta+1/{2}}} \eqsim \|u_0\|_{\DotH^{2\beta+1,p}/\bC}. \]
\end{cor}

Next, consider the inhomogeneous Cauchy problem ($u_0=0$, $F=0$).
\begin{theorem}[Well-posedness of the inhomogeneous Cauchy problem]
    \label{thm:wp-ic}
    Let $\beta>-1/2$ and $p_A(\beta)<p \le \infty$. Let $f \in T^p_\beta$. Then there exists a unique global weak solution $u$ to the \emph{inhomogeneous Cauchy problem}
    \begin{equation*}
        \tag{IC}
        \begin{cases}
            \partial_t u - \Div_x(A(t,x)\nabla_x u) = f, \quad (t,x) \in (0,\infty) \times \bR^n \\
            u(0)=0
        \end{cases}
    \end{equation*}
    with $\nabla u \in T^p_{\beta+1/2}$. Moreover, $u$ belongs to $C([0,\infty);\scrS') \cap T^p_{\beta+1}$ with
    \[ \|u\|_{T^p_{\beta+1}} + \|\nabla u\|_{T^p_{\beta+1/2}} \lesssim \|f\|_{T^p_\beta}. \]
\end{theorem} 

It is worth noting that our approach differs from the $L^p$-maximal regularity theory, which targets on finding \emph{mild} solutions to \eqref{e:ic} in a Banach space $X$ (\textit{e.g.}, $X=L^q(\bR^n)$) with estimates of $\partial_t u$ and $\Div(A\nabla u)$ in $L^p((0,T);X)$. The primary reason for this distinction is that $\DotH^{s,p}$ (and even $L^p$) spaces are not trace spaces in the sense of real interpolation theory. Besides, the characterization of $L^p$-maximal regularity for $A$ remains widely open. Existing affirmative results require additional Sobolev regularity assumptions, even in the simplest case $X=L^2(\bR^n)$ and $p=2$, see \cite{Dier-Zacher2017-MR-L2x-SobolevCoeff,Fackler2018-MR-Lpt-SobolevCoeff,Achache-Ouhabaz2019-MR-L2-H1/2Coeff}. The reader is also referred to \cite[\S 17.5]{Hytonen-vanNeerven-Veraar-Weis2023_BanachIII} for a more comprehensive survey.

At this stage, we do not have estimates for the tent space norms of $\partial_t u$ and $\Div(A\nabla u)$. However, for time-independent coefficients, it has been shown in \cite{Auscher-Hou2025_WPMRTent} that both $\partial_t u$ and $\Div(A(x)\nabla u)$ lie in $T^p_\beta$ with
\[ \|\partial_t u\|_{T^p_\beta} + \|\Div(A(x)\nabla u)\|_{T^p_\beta} \lesssim \|f\|_{T^p_\beta}. \]
See Remark \ref{rem:ic-real-max-reg} for detailed discussion.

The following theorem is concerned with the \emph{Lions equation} ($u_0=0$, $f=0$).
\begin{theorem}[Well-posedness of the Lions equation]
    \label{thm:wp-lions}
    Let $\beta>\beta_A$ and
    \[ \begin{cases}
        p^-_{\beta_A}(\beta) < p \le p^+_{\beta_A}(\beta) = \infty & \text{ if } \beta \ge -1/2 \\
        p^-_{\beta_A}(\beta) < p < p^+_{\beta_A}(\beta) & \text{ if } \beta_A < \beta < -1/2 
    \end{cases}. \]
    Let $F \in T^p_{\beta+1/2}$. Then there exists a unique global weak solution $u$ to the \emph{Lions equation}
    \begin{equation*}
        \tag{L}
        \begin{cases}
            \partial_t u - \Div_x(A(t,x)\nabla_x u) = \Div_x F, \quad (t,x) \in (0,\infty) \times \bR^n \\
            u(0)=0
        \end{cases}
    \end{equation*}
    with $\nabla u \in T^p_{\beta+1/2}$. Moreover, $u$ belongs to $C([0,\infty);\scrS') \cap T^p_{\beta+1}$ with
    \[ \|u\|_{T^p_{\beta+1}} + \|\nabla u\|_{T^p_{\beta+1/2}} \lesssim \|F\|_{T^p_{\beta+1/2}}. \]
\end{theorem}

To the best of our knowledge, for the Lions equation, only the case $\beta=-1/2$ has been addressed in a very recent work \cite{Auscher-Portal2025-Lions}, and only for $2 \le p \le \infty$. The case $p<2$ considered here appears to be new, since \cite{Auscher-Portal2025-Lions} only treats this case for \emph{real} coefficients and their argument relies on this condition in a crucial way.

Let us now present well-posedness for parabolic Cauchy problems of type \eqref{e:NaPC}, as announced.
\begin{theorem}[Well-posedness of non-autonomous parabolic Cauchy problems of divergence type]
    \label{thm:wp-NaPC}
    Let $\beta>\beta_A$ and
    \[ \begin{cases}
        p^-_{\beta_A}(\beta) < p \le p^+_{\beta_A}(\beta) = \infty & \text{ if } \beta \ge -1/2 \\
        p^-_{\beta_A}(\beta) < p < p^+_{\beta_A}(\beta) & \text{ if } \beta_A < \beta < -1/2 
    \end{cases}. \]
    Let $\gamma>-1/2$ and $p_A(\gamma)<q \le \infty$. Suppose
    \begin{equation}
        \label{e:wp-NaPC-embed-condition}
        \gamma \ge \beta, \quad 2\gamma-\frac{n}{q} = 2\beta-\frac{n}{p}.
    \end{equation}
    \begin{enumerate}[label=\normalfont(\roman*)]
        \item If $\beta<0$, then for any $u_0 \in \DotH^{2\beta+1,p}$, $F \in T^p_{\beta+1/2}$, and $f \in T^q_\gamma$, there exists a unique global weak solution $u$ to the Cauchy problem
        \[ \begin{cases}
            \partial_t u - \Div_x( A(t,x) \nabla_x u ) = f + \Div_x F, \quad (t,x) \in (0,\infty) \times \bR^n \\
            u(0)=u_0
        \end{cases} \]
        so that $\nabla u \in T^p_{\beta+1/2}$. Moreover, the estimate holds that
        \[ \|\nabla u\|_{T^p_{\beta+1/2}} \lesssim \|u_0\|_{\DotH^{2\beta+1,p}} + \|F\|_{T^p_{\beta+1/2}} + \|f\|_{T^q_\gamma}, \]
        and $u$ belongs to $C([0,\infty);\scrS')$. If $u_0=0$, then $u$ also belongs to $T^p_{\beta+1}$ with
        \[ \|u\|_{T^p_{\beta+1}} \lesssim \|F\|_{T^p_{\beta+1/2}} + \|f\|_{T^q_\gamma}. \]

        \item If $\beta \ge 0$, then the above statements also hold when $u_0$ is constant (for which $\|u_0\|_{\DotH^{2\beta+1,p}}=0$).
    \end{enumerate}
\end{theorem}
The proof is a straightforward combination of Theorems \ref{thm:wp-hc}, \ref{thm:wp-ic}, and \ref{thm:wp-lions}. We only need to mention that the condition \eqref{e:wp-NaPC-embed-condition} allows one to use Sobolev embedding of tent spaces (see Section \ref{ssec:tent}) to show that the weak solution $u_1$ to the inhomogeneous Cauchy problem \eqref{e:ic} with $f \in T^q_{\gamma}$ satisfies $\nabla u_1 \in T^p_{\beta+1/2}$ and
\[ \|\nabla u_1\|_{T^p_{\beta+1/2}} \lesssim \|\nabla u_1\|_{T^q_{\gamma+1/2}} \lesssim \|f\|_{T^q_\gamma}. \]

Notice that Theorem \ref{thm:rep} \ref{item:rep_beta>0} implies for $\beta \ge 0$, the only initial data compatible with the solution class $\nabla u \in T^p_{\beta+1/2}$ are constant.

To finish the introduction, let us mention an interesting corollary that we did not find in the literature. Note that for $p=2$, the tent space $T^2_\beta$ identifies with the time-weighted $L^2$-space
\[ T^2_\beta \simeq L^2\left( (0,\infty),t^{-2\beta} dt;~L^2(\bR^n) \right) =: L^2_\beta(\bR^{1+n}_+), \quad \beta \in \bR. \]
In this special case, we have
\begin{cor}
    Let $\beta>-1/2$. For any $F \in L^2_{\beta+1/2}(\bR^{1+n}_+)$ and $f \in L^2_\beta(\bR^{1+n}_+)$, there exists a unique global weak solution $u$ to the Cauchy problem
    \[ \begin{cases}
        \partial_t u - \Div_x( A(t,x) \nabla_x u ) = f + \Div_x F, \quad (t,x) \in (0,\infty) \times \bR^n \\
        u(0)=0
    \end{cases} \]
    so that $\nabla u \in L^2_{\beta+1/2}(\bR^{1+n}_+)$. Moreover, $u$ lies in $L^2_{\beta+1}(\bR^{1+n}_+)$ with
    \[ \|u\|_{L^2_{\beta+1}(\bR^{1+n}_+)} + \|\nabla u\|_{L^2_{\beta+1/2}(\bR^{1+n}_+)} \lesssim \|F\|_{L^2_{\beta+1/2}(\bR^{1+n}_+)} + \|f\|_{L^2_{\beta}(\bR^{1+n}_+)}. \]
    The same results are also valid for $\beta_A<\beta \le -1/2$ and $f=0$.
\end{cor}

\subsection{Organization}
\label{ssec:organization}
The paper is organized as follows.

Section \ref{sec:spaces} collects definitions and fundamental properties of function spaces to be used. 

Section \ref{sec:fact} recalls the $L^2$-theory of weak solutions and propagators. We also establish the $L^p$-theory of propagators in Proposition \ref{prop:Gamma-Lp}, using the notion of $L^p-L^q$ off-diagonal decay (see Definition \ref{def:odd}).

Sections \ref{sec:ic}, \ref{sec:lions}, and \ref{sec:hc} are devoted to proving existence of weak solutions to \eqref{e:ic}, \eqref{e:lions}, and \eqref{e:hc}, as announced in Theorems \ref{thm:wp-ic}, \ref{thm:wp-lions}, and \ref{thm:wp-hc}, respectively. Main results of these three sections are Theorems \ref{thm:ic} (inhomogeneous Cauchy problems), \ref{thm:lions} (the Lions equation), and \ref{thm:hc} (homogeneous Cauchy problems). Section \ref{sec:unique} is concerned with uniqueness and representation. The proof of Theorem \ref{thm:rep} is presented. 

Section \ref{sec:Besov} contains extension to homogeneous Besov spaces.

\subsection{Notation}
\label{ssec:notation}
Throughout the paper, for any $q,r \in (0,\infty]$, we write 
\[ [q,r]:=\frac{1}{q}-\frac{1}{r}, \]
if there is no confusion with closed intervals. We say $X \lesssim Y$ (or $X \lesssim_A Y$, resp.) if $X \le CY$ with an irrelevant constant $C$ (or depending on $A$, resp.), and say $X \eqsim Y$ if $X \lesssim Y$ and $Y \lesssim X$. 

Write $\bR^{1+n}_+:= (0,\infty) \times \bR^n$. For any (Euclidean) ball $B \subset \bR^n$, write $r(B)$ for the radius of $B$. For any function $f$ defined on $\bR^{1+n}_+$, denote by $f(t)$ the function $x \mapsto f(t,x)$ for any $t>0$.

Let $(X,\mu)$ be a measure space. For any measurable subset $E \subset X$ with finite positive measure and $f \in L^1(E,\mu)$, we write
\[ \fint_E f d\mu := \frac{1}{\mu(E)} \int_E f d\mu. \]

Often, we omit the unweighted Lebesgue measure in the integral and the domain in the function space, if it is clear from the context. We use the sans-serif font $\ssfc$ in the scripts of function spaces in short of ``with compact support" in the prescribed set, and ${\loc}$ if the prescribed property holds on all compact subsets of the prescribed set.

\subsection{Acknowledgments}
\label{ssec:acknowledg}
The author would like to thank his Ph.D. advisor Pascal Auscher for enlightening discussions and generous help.
\section{Function spaces}
\label{sec:spaces}

\subsection{Homogeneous Hardy--Sobolev spaces}
\label{ssec:Hsp-Epdelta}
Readers can refer to \cite[\S 5]{Triebel1983-book} for definitions and basic properties of Triebel--Lizorkin spaces and homogeneous Hardy--Sobolev spaces, as well as the proof of the following facts. Some needed refinements can be found in \cite[\S 2.4.3 \& 3.3.3]{Sawano2018_Besov} and the references therein, see also \cite[\S 2.1]{Auscher-Hou2024-HCL}. 

Denote by $\scrS(\bR^n)$ the space of Schwartz functions and by $\scrS'(\bR^n)$ the space of tempered distributions. For any $f \in \scrS'(\bR^n)$, write $\cF(f)$ for the Fourier transform of $f$. 

Denote by $C$ the annulus $\{\xi \in \bR^n:2^{-1} \le |\xi| \le 2^2\}$ and by $2^j C$ the annulus $\{\xi \in \bR^n:2^{j-1} \le |\xi| \le 2^{j+2}\}$ for any $j \in \bZ$. Pick $\chi \in \Cc(\bR^n)$ with $\supp(\chi) \subset C$ and $\sum_{j \in \bZ} \chi(2^{-j} \xi)=1$ for any $\xi \ne 0$. Let $\Delta_j$ be the $j$-th \textit{Littlewood--Paley operator (associated with $\chi$)} given by
\[ \Delta_j f := \cF^{-1} (\chi(2^{-j} \cdot) \cF(f)), \quad f \in \scrS'. \]

Let $\scrP$ be the space of polynomials $\bC[x_1,\dots,x_n]$, $\scrP_0:=\{0\}$, and $\scrP_m$ be the subspace of $\scrP$ consisting of polynomials of degree less than $m$ for $m \ge 1$.

Let $s \in \bR$ and $0<p<\infty$. The \emph{Triebel--Lizorkin space} $\DotF^s_{p,2}$ consists of $f \in \scrS'/\scrP$ for which
\[ \|f\|_{\DotF^s_{p,2}} := \left \| \left ( \sum_j |2^{js} \Delta_j f|^2 \right )^{1/2} \right \|_{L^p} < \infty. \]
For $p=\infty$, it is defined similarly by approximation. 

Define $\nu(s,p):=\max\{ 0,[ s-\frac{n}{p} ]+1 \}$. For any $f \in \DotF^s_{p,2}$, the Littlewood--Paley series $\sum_j \Delta_j f$ converges in $\scrS'/\scrP_{\nu(s,p)}$. It induces an isometric embedding $\iota: \DotF^s_{p,2} \to \scrS'/\scrP_{\nu(s,p)}$ given by
\[ \iota(f):=\sum_j \Delta_j f. \]

\begin{definition}[Homogeneous Hardy--Sobolev spaces]
    \label{def:Hsp}
    Let $s \in \bR$ and $0<p \le \infty$. The \emph{homogeneous Hardy--Sobolev space} $\DotH^{s,p}$ consists of $f \in \scrS'$ whose class in $\scrS'/\scrP_{\nu(s,p)}$ belongs to $\iota(\DotF^s_{p,2})$. The (quasi-semi-)\ norm of $\DotH^{s,p}$ is given by
    \[ \|f\|_{\DotH^{s,p}} := \|[f]\|_{\DotF^s_{p,2}}, \]
    where $[f]$ denotes the class of $f$ in $\scrS'/\scrP$.
\end{definition}

For $s=0$, $\DotH^{0,p}$ coincides with the Hardy space $H^p$ if $p \le 1$, the Lebesgue space $L^p$ if $1<p<\infty$, and $\bmo$ if $p=\infty$, up to equivalent (quasi-/semi-)norms. Moreover, for any $s \in \bR$, $\DotH^{s,\infty}/\scrP_{\nu(s,\infty)}$ is isomorphic to $\bmo^s$ introduced by Strichartz \cite{Strichartz1980_BMOs}. 

Our homogeneous Hardy--Sobolev spaces $\DotH^{s,p}$ can be regarded as ``realizations'' of Triebel--Lizorkin spaces $\DotF^s_{p,2}$ in $\scrS'$, and hence inherit many fundamental properties, such as duality, complex interpolation, and Sobolev embedding. In particular, let us mention that for any $s \in \bR$, $\DotH^{s,p} \cap L^2$ is dense in $\DotH^{s,p}$ if $p<\infty$, or weak*-dense if $p=\infty$.

Let $\DotH^{s,p}+\bC$ be the collection of $f \in \scrS'$ of the form $f=g+c$ for some $g \in \DotH^{s,p}$ and $c \in \bC$, endowed with the (quasi-semi-)norm
\[ \|f\|_{\DotH^{s,p}/\bC} := \inf_g \|g\|_{\DotH^{s,p}}, \]
where the infimum is taken among all $g \in \DotH^{s,p}$ for which $f=g+c$ for some $c \in \bC$. Note that when $\nu(s,p) \ge 0$, $\DotH^{s,p}$ contains constants and $\|f+c\|_{\DotH^{s,p}} = \|f\|_{\DotH^{s,p}}$ for any $c \in \bC$, so $\DotH^{s,p}+\bC = \DotH^{s,p}$.

\subsection{Tent spaces}
\label{ssec:tent}
We adapt the original definitions in \cite{Coifman-Meyer-Stein1985_TentSpaces} to the parabolic settings. Let $\beta \in \bR$ and $0<p<\infty$. The \emph{(parabolic) tent space} $T^p_\beta$ consists of measurable functions $f$ on $\bR^{1+n}_+$ for which
\[ \|f\|_{T^{p}_\beta} := \left ( \int_{\bR^n} \left ( \int_0^\infty \fint_{B(x,t^{1/2})} |t^{-\beta} f(t,y)|^2 \, dtdy \right )^{p/2} dx \right )^{1/p} < \infty. \]
For $p=\infty$, the space $T^{\infty}_{\beta}$ consists of measurable functions $f$ on $\bR^{1+n}_+$ for which
\[ \|f\|_{T^{\infty}_{\beta}} := \sup_{B} \left ( \int_0^{r(B)^2} \fint_B |t^{-\beta} f(t,y)|^2\,  dtdy \right )^{1/2} < \infty,  \]
where $B$ describes balls in $\bR^n$. We list below some properties to be used. The reader can refer to \cite{Coifman-Meyer-Stein1985_TentSpaces,Hofmann-Mayboroda-McIntosh2011-HpL,Amenta-Auscher2018-book,Amenta2018_WeightedTent,Auscher-Hou2025_WPMRTent} for the proofs.
\begin{enumerate}[label=\normalfont(\roman*)]
    \item Let $\beta \in \bR$ and $0<p \le \infty$. Then $T^p_\beta$ embeds into $L^2_{\loc}(\bR^{1+n}_+)$, and $L^2_\ssfc(\bR^{1+n}_+)$ is dense (resp. weak*-dense) in $T^p_\beta$ if $0<p<\infty$ (resp. $p=\infty$).

    \item (Duality) Let $\beta \in \bR$ and $1 \le p < \infty$. The dual of $T^p_\beta$  identifies with $T^{p'}_{-\beta}$ via $L^2(\bR^{1+n}_+)$-duality.

    \item (Interpolation) Let $\beta_0, \beta_1 \in \bR$, $0<p_0<p_1 \le \infty$. Suppose $1/p=(1-\theta)/p_0+\theta/p_1$ and $\beta=(1-\theta)\beta_{0}+\theta\beta_{1}$ for some $\theta \in (0,1)$. Then the complex interpolation space $[T^{p_0}_{\beta_{0}}, T^{p_1}_{\beta_{1}}]_\theta$ identifies with $T^p_\beta$.

    \item (Embedding) Let $\beta_0>\beta_1$ and $0<p_0<p_1 \le \infty$. Suppose
    \[ 2\beta_0 - \frac{n}{p_0} = 2\beta_1 - \frac{n}{p_1}. \]
    Then $T^{p_0}_{\beta_0}$ embeds into $T^{p_1}_{\beta_1}$.
\end{enumerate}

Let us also introduce two scales of retractions of tent spaces.

\subsubsection{Slice spaces}
\label{sssec:slice}
Slice spaces are introduced by \cite[\S 3]{Auscher-Mourgoglou2019_Ep_delta}. Let $0< p \le \infty$ and $\delta>0$. The \emph{(parabolic) slice space} $E^p_\delta$ consists of measurable functions $f$ on $\bR^n$ for which
\[ \|f\|_{E^p_\delta} := \left\| \left( \fint_{B(\cdot,\delta^{1/2})} |f(y)|^2 dy \right)^{1/2} \right\|_{L^p} < \infty. \]
For any $\delta>0$, $E^p_\delta$ is a retraction of $T^p_0$. We mention here some other fundamental properties.
\begin{enumerate}[label=(\roman*)]
    \item (Duality) For $1 \le p < \infty$, the dual of $E^p_\delta$ identifies with $E^{p'}_\delta$ via $L^2(\bR^n)$-duality.

    \item (Interpolation) Let $0<p_0<p_1 \le \infty$. Suppose $1/p=(1-\theta)/p_0+\theta/p_1$ for some $\theta \in (0,1)$. Then the complex interpolation space $[E^{p_0}_\delta, E^{p_1}_\delta]_\theta$ identifies with $E^p_\delta$. 

    \item (Change of aperture) Let $\delta'>0$. Then $E^p_\delta$ agrees with $E^p_{\delta'}$ as sets, with equivalents norms. More precisely,
    \begin{equation}
        \label{e:Ep-aperture}
        \min\left\{ 1,\left( \frac{\delta}{\delta'} \right)^{ \frac{n}{2} (\frac{1}{2}-\frac{1}{p}) } \right\} \|f\|_{E^p_{\delta}} \lesssim_{n,p} \|f\|_{E^p_{\delta'}} \lesssim_{n,p} \max\left\{ 1,\left( \frac{\delta}{\delta'} \right)^{ \frac{n}{2} (\frac{1}{2}-\frac{1}{p}) } \right\} \|f\|_{E^p_{\delta}}.
    \end{equation}
\end{enumerate}

\subsubsection{Divergence of tent spaces}
\label{sssec:divTpbeta}
Let $\beta \in \bR$ and $1<p<\infty$. Let $\Div T^p_{\beta}$ be the collection of distributions $f \in \scrD'(\bR^{1+n}_+)$ of the form $f=\Div F$ for some $F \in T^p_{\beta}$, endowed with the norm
\[ \|f\|_{\Div T^p_{\beta}} := \inf \left\{ \|F\|_{T^p_{\beta}} : \Div F = f, ~ F \in T^p_{\beta} \right\}. \] 

Denote by $R$ the Riesz transform on $L^2(\bR^{1+n}_+)$ given by
\[ R(G)(t) = \nabla (-\Delta)^{-1/2} G(t), \quad t>0. \]
We infer from \cite[Theorem 2.4]{Auscher-Prisuelos-Arribas2017-extrapolation} that $R$ is bounded on $T^p_{\beta}$. By duality, its dual operator $R^*$, given by
\[ R^*(G)(t) := -(-\Delta)^{-1/2} \Div G(t), \quad t>0, \]
is also bounded on $T^p_{\beta}$. Let $f \in \Div T^p_{\beta}$. Pick $F \in T^p_{\beta}$ with $\Div F = f$. Define the operator $S:\Div T^p_{\beta} \to T^p_{\beta}$ by
\[ S(f) := RR^*(F) = -R(-\Delta)^{-1/2} \Div (F). \]
The definition is independent of the choice of $F$. Using boundedness of $R$ and $R^\ast$, we get $\|S(f)\|_{T^p_{\beta}} \lesssim \|F\|_{T^p_{\beta}}$. Taking infimum among all such $F$ yields $\|S(f)\|_{T^p_{\beta}} \lesssim \|f\|_{\Div T^p_{\beta}}$, so $S$ is bounded.

\begin{lemma}
    \label{lem:divT}
    Let $\beta \in \bR$ and $1<p<\infty$. Then $\Div S = \id_{\Div T^p_{\beta}}$. Consequently, $\Div T^p_{\beta}$ is a Banach space, and $\Div: T^p_{\beta} \to \Div T^p_{\beta}$ is a retraction with $S:\Div T^p_{\beta} \to T^p_{\beta}$ as the coretraction.
\end{lemma}

\begin{proof}
    We only need to prove the identity. Let $f \in \Div T^p_{\beta}$ and $F \in T^p_{\beta}$ with $\Div F = f$. Pick $F_n \in L^2_\ssfc((0,\infty);W^{1,2}_\ssfc(\bR^n))$ as a sequence approximating $F$ in $T^p_{\beta}$. For any $t>0$, we have $\Div S (\Div F_n)(t) = -\Div \nabla (-\Delta)^{-1/2} (-\Delta)^{-1/2} \Div F_n(t) = \Div F_n(t)$ in $L^2(\bR^n)$, so 
    \[ \Div S(\Div F_n) = \Div F_n \quad \text{ in } L^2(\bR^{1+n}_+). \]
    Ensured by boundedness of $\Div$ and $S$, we take limits on both sides to get $\Div S(f) = \Div S(\Div F) = \Div F = f$ in $\Div T^p_{\beta}$.
\end{proof}
\section{Basic facts of weak solutions}
\label{sec:fact}

We first recall the definition of weak solutions. Denote by $W^{1,2}$ the \emph{(inhomogeneous) Sobolev space} endowed with the norm $\|f\|_{W^{1,2}}:=\|f\|_{L^2}+\|\nabla f\|_{L^2}$.
\begin{definition}[Weak solutions]
    \label{def:weak-sol}
    Let $\Omega$ be an open subset of $\bR^n$ and $0 \le a < b \le \infty$. Let $f$ and $F$ be in $\scrD'((a,b) \times \Omega)$. A function $u \in L^2_{\loc}((a,b);W^{1,2}_{\loc}(\Omega))$ is called a \emph{weak solution} to the equation
    \[ \partial_t u - \Div (A\nabla u) = f + \Div F \]
    with \emph{source term} $f+\Div F$, if for any $\phi \in \Cc((a,b) \times \Omega)$,
    \begin{equation}
        \label{e:weak-sol}
        -\int_{(a,b) \times \Omega} u ~ \partial_t \phi + \int_{(a,b) \times \Omega} (A\nabla u) \cdot \nabla \phi = (f,\phi) - (F,\nabla \phi).
    \end{equation}
    The pairs on the right-hand side stand for pairing of distributions and test functions on $(a,b) \times \Omega$. We say $u$ is a \emph{global} weak solution if \eqref{e:weak-sol} holds for $(a,b) \times \Omega = (0,\infty) \times \bR^n$.

    Let $u_0 \in \scrD'(\Omega)$. We say $u$ satisfies the \emph{initial condition} $u(a)=u_0$ if $u(t)$ converges to $u_0$ in $\scrD'(\Omega)$ as $t \to a+$.
\end{definition}
There is a corresponding definition of (global) weak solutions to the backward equation $-\partial_s u - \Div(A(s)^\ast \nabla u) = f + \Div F$. We leave the precise formulation to the reader.

\subsection{Energy inequality}
\label{ssec:Caccioppoli}
We recall without proof a form of Caccioppoli's inequality.
\begin{lemma}[Caccioppoli's inequality]
    \label{lem:Caccioppoli}
    Let $0<a<b<\infty$ and $B \subset \bR^n$ be a ball. Let $f$ and $F$ be in $L^2((a,b) \times 2B)$. Let $u \in L^2((a,b);W^{1,2}(2B))$ be a weak solution to $\partial_t u - \Div(A\nabla u) = f + \Div F$ in $(a,b) \times 2B$. Then $u$ belongs to $C([a,b];L^2(B))$ with
    \begin{align*}
        \|u(b)\|_{L^2(B)}^2 \lesssim
        &\left ( \frac{1}{r(B)^2}+\frac{1}{b-a} \right ) \int_a^b \|u(s)\|_{L^2(2B)}^2\, ds \\
        &+ r(B)^2 \int_a^b \|f(s)\|_{L^2(2B)}^2 \,ds + \int_a^b \|F(s)\|_{L^2(2B)}^2\, ds .
    \end{align*}
    Moreover, for any $c \in (a,b)$, it holds that
    \begin{align*}
        \int_c^b \|\nabla u(s)\|_{L^2(B)}^2 \, ds \lesssim &\frac{1}{c-a} \left ( 1+\frac{b-a}{r(B)^2} \right ) \int_a^b \|u(s)\|_{L^2(2B)}^2 \, ds \\
        &+ \frac{r(B)^2 (b-a)}{c-a} \int_a^b \|f(s)\|_{L^2(2B)}^2 \,ds \\
        &+ \frac{b-a}{c-a} \int_a^b \|F(s)\|_{L^2(2B)}^2 \, ds.
    \end{align*}
\end{lemma}
There is a corresponding version for weak solutions to the backward equation $-\partial_s u - \Div(A(s)^\ast \nabla u) = f + \Div F$. We refer to ``Caccioppoli's inequality" in both cases in the sequel. 
\begin{cor}[\textit{A priori} tent space estimates]
    \label{cor:Caccioppoli-tent}
    Let $\beta \in \bR$ and $0<p \le \infty$. Let $u$ be a global weak solution to $\partial_t u - \Div(A \nabla u) = f + \Div F$. Then the \textit{a priori} estimate holds that
    \[ \|\nabla u\|_{T^p_{\beta+1/2}} \lesssim \|u\|_{T^p_{\beta+1}} + \|F\|_{T^p_{\beta+1/2}} + \|f\|_{T^p_{\beta}}. \]
    This inequality also holds for weak solutions to the backward equation.
\end{cor}
\begin{proof}
    The proof follows from the same arguments for time-independent coefficients, by applying Caccioppoli's inequality to averages on local Whitney cubes $(t,2t) \times B(x,t^{1/2})$, see \cite[Corollary 4.3]{Auscher-Hou2024-HCL}.
\end{proof}

\subsection{Propagators}
\label{ssec:propagator}
We first recall the $L^2$-theory.

\begin{prop}[$L^2$-theory]
    \label{prop:L2}
    Let $u_0 \in L^2$ and $F \in L^2(\bR^{1+n}_+) \simeq T^2_0$. Then there exists a unique global weak solution $u$ to the Cauchy problem
    \[ \begin{cases}
        \partial_t u - \Div(A\nabla u) = \Div F \\
        u(0) = u_0
    \end{cases} \]
    so that $\nabla u \in L^2(\bR^{1+n}_+)$. Moreover, $u$ belongs to $C([0,\infty);L^2)$ with
    \[ \sup_{t \ge 0} \|u(t)\|_{L^2} \le \|u_0\|_{L^2} + (2\Lambda_0)^{-1/2} \|F\|_{L^2(\bR^{1+n}_+)}, \]
    and
    \[ \|\nabla u\|_{L^2(\bR^{1+n}_+)} \le (2\Lambda_0)^{-1/2} \|u_0\|_{L^2} + \Lambda_0^{-1} \|F\|_{L^2(\bR^{1+n}_+)}. \]
\end{prop}
The existence is due to \cite{Lions1957_L2}, and the uniqueness in the class $\nabla u \in L^2(\bR^{1+n}_+)$ is due to \cite{Auscher-Monniaux-Portal2019Lp}, see \cite[Theorem 2.2]{Auscher-Portal2025-Lions} for a survey. Moreover, \cite{Auscher-Monniaux-Portal2019Lp} further exploits the $L^2$-theory to establish the notion of propagators. The reader can refer there for the proof of the following two corollaries.

\begin{cor}[Propagators]
    \label{cor:Gamma}
    There exists a family of contractions on $L^2(\bR^n)$, $(\Gamma_A(t,s))_{0 \le s \le t < \infty}$, called \emph{(forward) propagators associated to $A$} so that for any $h \in L^2$ and $s \ge 0$, $u(t):=\Gamma_A(t,s) h$ is the unique weak solution on $(s,\infty) \times \bR^n$ to the Cauchy problem
    \[ \begin{cases}
        \partial_t u - \Div(A\nabla u) = 0, \quad t>s \\
        u(s)=h
    \end{cases} \]
    with $\nabla u \in L^2((s,\infty) \times \bR^n)$. Moreover, $u$ belongs to $C_0([s,\infty);L^2)$.
    \footnote{Here, $C_0([s,\infty);E)$ is the space of continuous functions with limit 0 as $t \to \infty$ in the prescribed topology.}
\end{cor} 

By reversing time, we also obtain the notion of propagators for the backward equation, called backward propagators.
\begin{cor}[Backward propagators]
    \label{cor:Gamma-}
    There exists a family of contractions on $L^2$, $(\Gamma_{A^\ast}^-(s,t))_{0 \le s \le t < \infty}$, called \emph{backward propagators associated to $A^\ast$}, so that for any $h \in L^2$ and $t>0$, $\tilu(s):=\Gamma_{A^\ast}^-(s,t) h$ is the unique weak solution on $(0,t) \times \bR^n$ to the backward Cauchy problem
    \[ \begin{cases}
        -\partial_s \tilu - \Div(A(s)^\ast \nabla \tilu) = 0, \quad 0<s<t \\
        u(t)=h
    \end{cases} \]
    so that $\nabla \tilu \in L^2((0,t) \times \bR^n)$. Moreover, for fixed $t>0$, it satisfies
    \begin{equation}
        \label{e:Gamma-dual}
        \Gamma^-_{A^\ast}(s,t) = \Gamma_{A}(t,s)^\ast = \Gamma_{\tilA_t}(t-s,0), \quad 0 \le s \le t,
    \end{equation}
    where
    \begin{equation}
        \label{e:tilAt}
        \tilA_t(s,x) := 
        \begin{cases}
            A^\ast(t-s,x), & \text{ if } 0 \le s \le t \\
            \Lambda_0 \bI & \text{ if } s>t
        \end{cases}.
    \end{equation}
    Consequently, $\tilu$ belongs to $C([0,t];L^2)$.
\end{cor}
Remark that $\tilA_t$ has the same ellipticity as $A$. When $s>t$, our construction of $\tilA_t$ is different from the original proof in \cite[Lemma 3.16]{Auscher-Monniaux-Portal2019Lp}, but it is irrelevant for the identity \eqref{e:Gamma-dual}, which holds for $0 \le s \le t$. 

The propagators provide explicit formulas for weak solutions. Recall the propagator solution map $\cE_A$ defined in \eqref{e:EA} from $L^2$ to $L^\infty((0,\infty);L^2) \cap L^2_{\loc}((0,\infty);W^{1,2}_{\loc})$ by
\[ \cE_A(u_0)(t,x) := (\Gamma_A(t,0)u_0)(x). \]

Define the \emph{Duhamel operator} $\cL^A_1$ from $L^2(\bR^{1+n}_+)$ to $L^\infty_{\loc}([0,\infty);L^2)$ by the $L^2$-valued Bochner integrals (verified below)
\begin{equation}
    \label{e:LA1}
    \cL^A_1(f)(t) := \int_0^t \Gamma_A(t,s) f(s) ds, \quad t>0,
\end{equation}
and the \emph{backward Duhamel operator} $(\cL^A_1)^\ast$ from $L^2_{\ssfc}(\bR^{1+n}_+)$ to $L^\infty_{\loc}([0,\infty);L^2)$ by the $L^2$-valued Bochner integrals (also verified below)
\begin{equation}
    \label{e:LA1*}
    (\cL^A_1)^\ast(f)(s) := \int_s^\infty \Gamma_A(t,s)^\ast f(t) dt, \quad 0<s<t.
\end{equation}
As we shall see in Lemma \ref{lem:LA-SIO}, $(\cL^A_1)^\ast$ is indeed the adjoint of $\cL^A_1$ with respect to the $L^2(\bR^{1+n}_+)$-duality.

Let us explain why the integral in \eqref{e:LA1} is a Bochner integral. The one in \eqref{e:LA1*} follows similarly. We first verify the strong measurability of the function $s \mapsto \Gamma_A(t,s) f(s)$, valued in $L^2$. Indeed, we infer from Corollary \ref{cor:Gamma-} that for any $\phi,\psi \in L^2$, the function
\begin{equation}
    \label{e:s->Gamma(t,s)-weak-meas}
    s \mapsto \langle \Gamma_A(t,s) \phi, \psi \rangle = \langle \phi, \Gamma_{A}(s,t)^\ast \psi \rangle = \langle \phi, \Gamma^-_{A^\ast}(s,t) \psi \rangle
\end{equation}
is continuous, hence (Borel) measurable. So for any $f \in L^2(\bR^{1+n}_+)$, $s \mapsto \langle \Gamma_A(t,s) f(s),\psi \rangle$ is measurable. We thus get $s \mapsto \Gamma_A(t,s)f(s)$ is weakly measurable, and hence strongly measurable by Pettis' measurability theorem (see \cite[Theorem 1.1.20]{Hytonen-vanNeerven-Veraar-Weis2016-Banach-I}), since $L^2$ is separable. The integrability comes from the fact that $\|\Gamma_A(t,s)\|_{\scrL(L^2)} \le 1$.

Define the \emph{Lions operator} $\cR^A_{1/2}$ by the formal integral
\begin{equation}
    \label{e:RA1/2}
    \cR^A_{1/2}(F)(t) := \int_0^t \Gamma_A(t,s) \Div F(s) ds, \quad t>0.
\end{equation}
We do not know whether the integral in \eqref{e:RA1/2} converges in the sense of Bochner integrals on $L^2$, because we do not have estimates for the operator norm of  $\Gamma_A(t,s)\Div$ in $\scrL(L^2)$ for each $s$ and $t$. Instead, it is only interpreted in the weak sense, \textit{i.e.}, for any $t>0$, $\cR^A_{1/2}(F)(t)$ is defined as a continuous linear functional on $L^2$ by
\[ \langle \cR^A_{1/2}(F)(t),h \rangle_{L^2(\bR^n)} := -\int_0^t \langle F(s), \nabla \Gamma_A(t,s)^\ast h \rangle_{L^2(\bR^n)} ds, \]
for any $h \in L^2$. The integral on the right-hand side converges as
\begin{align*}
    \int_0^t \left| \langle F(s), \nabla \Gamma_A(t,s)^\ast h \rangle_{L^2(\bR^n)} \right| ds 
    &\le \|F\|_{L^2(\bR^{1+n}_+)} \|\nabla \Gamma_A(t,\cdot)^\ast h\|_{L^2((0,t) \times \bR^n)} \\
    &\le (2\Lambda_0)^{-1/2} \|F\|_{L^2(\bR^{1+n}_+)} \|h\|_{L^2}.
\end{align*}
In the last inequality, we use \eqref{e:Gamma-dual} and the energy inequality in Proposition \ref{prop:L2} for reversed time, noting that $\tilA_t$ has the same ellipticity as $A$. This estimate also yields that $\cR^A_{1/2}$ is bounded from $L^2(\bR^{1+n}_+)$ to $L^\infty((0,\infty);L^2)$.

\begin{prop}
    \label{prop:formula}
    Let $u_0 \in L^2$, $f \in L^2_\ssfc(\bR^{1+n}_+)$, and $F \in L^2(\bR^{1+n}_+)$. 

    \begin{enumerate}[label=\normalfont(\roman*)]
        \item \label{item:L2-ic}
        $v:=\cL^A_1(f)$ is a global weak solution to the inhomogeneous Cauchy problem
        \begin{equation}
            \tag{IC}
            \label{e:ic}
            \begin{cases}
                \partial_t v - \Div(A\nabla v) = f \\
                v(0)=0
            \end{cases}.
        \end{equation}

        \item \label{item:L2-ic-back}
        $w:=(\cL^A_1)^\ast(f)$ is a weak solution on $\bR^{1+n}_+$ to the backward equation
        \[ -\partial_s w - \Div(A(s)^\ast\nabla w) = f. \]

        \item \label{item:L2-HCL}
        Let $u$ be the unique global weak solution to the Cauchy problem
        \[ \begin{cases}
            \partial_t u - \Div(A\nabla u) = \Div F \\
            u(0) = u_0
        \end{cases} \]
        with $\nabla u \in L^2(\bR^{1+n}_+)$. Then the following properties hold.
        \begin{enumerate}[label=\normalfont(\arabic*)]
            \item \label{item:formula-L2-Duhamel}
            {\normalfont (Duhamel's formula)}
            $u = \cE_A(u_0) + \cR^A_{1/2}(F)$ in $\scrD'(\bR^{1+n}_+)$.
            
            \item \label{item:formula-L2-RA1/2}
            Define $\tilF:=(A-\bI) \nabla \cR^A_{1/2}(F) + F$ in $L^2(\bR^{1+n}_+)$. Then
            \[ \cR^A_{1/2}(F) = \cR^{\bI}_{1/2}(\tilF) = \Div \cL^{\bI}_1(\tilF) \quad \text{ in }  \scrD'(\bR^{1+n}_+). \]

            \item \label{item:formula_L2_EA}
            It holds that
            \begin{align*}
                \cE_A(u_0) 
                &= \cE_{\bI}(u_0) + \cR^A_{1/2}((A-\bI) \nabla \cE_{\bI}(u_0)) \\
                &= \cE_{\bI}(u_0) + \cR^{\bI}_{1/2}((A-\bI) \nabla \cE_A(u_0)).
            \end{align*}
        \end{enumerate}
    \end{enumerate}
\end{prop}

\begin{proof}
    The statement \ref{item:L2-ic} directly follows from \cite[Theorem 2.54]{Auscher-Egert2023-Propagator}. Applying the same argument to the backward equation gives \ref{item:L2-ic-back}, thanks to \eqref{e:Gamma-dual}. To prove \ref{item:L2-HCL} \ref{item:formula-L2-Duhamel}, we infer from \cite[Theorem 2.54]{Auscher-Egert2023-Propagator} that Duhamel's formula holds in $L^\infty((0,T);L^2)$ for any $T>0$, hence in $\scrD'(\bR^{1+n}_+)$. The statements \ref{item:formula-L2-RA1/2} and \ref{item:formula_L2_EA} follow from the same arguments as in \cite[Corollary 4.5]{Auscher-Hou2024-HCL} for time-independent coefficients. Details are left to the reader.
\end{proof}

Let us formulate the $L^p$-theory of the propagators by using the notion of off-diagonal decay introduced by \cite{Auscher-Hou2025_WPMRTent}. Denote by $\Delta$ the set $\{(t,s) \in (0,\infty) \times (0,\infty): t=s\}$ and write $\Delta^c := ( (0,\infty) \times (0,\infty) ) \setminus \Delta$.

\begin{definition}[Off-diagonal decay]
    \label{def:odd}
    Let $1 \le p \le q \le \infty$. Let $\{K(t,s)\}_{(t,s) \in \Delta^c}$ be a family of bounded operators on $L^2(\bR^n)$. We say $\{K(t,s)\}$ satisfies the \emph{(exponential) $L^p-L^q$ off-diagonal decay} if there are constants $c,C>0$ so that for any $t>s$, $E,F \subset \bR^n$ as Borel sets, and $f \in L^2 \cap L^p$,
    \[ \|\I_E K(t,s) \I_F f\|_{L^q} \le C (t-s)^{-\frac{n}{2}[p,q]} \exp\left( -c \frac{\dist(E,F)^2}{t-s} \right) \|\I_F f\|_{L^p}. \]
\end{definition}

\begin{prop}[$L^p$-theory of propagators]
    \label{prop:Gamma-Lp}
    Let $A \in L^\infty(\bR^{1+n}_+;\mat_n(\bC))$ be uniformly elliptic. Let $(\Gamma_A(t,s))$ be the propagators associated to $A$. Then the following properties hold.
    \begin{enumerate}[label=\normalfont(\roman*)]
        \item \label{item:p-+(A)-bd}
        There exist $p_-(A) \in [1,2)$ and $p_+(A) \in (2,\infty]$ so that $(p_-(A),p_+(A))$ is the maximal open set of exponents $p \in [1,\infty]$ for which $(\Gamma_A(t,s))$ is uniformly bounded on $L^p$.

        \item \label{item:p-+(A)-odd}
        For $p_-(A)<p \le 2$, $(\Gamma_A(t,s))$ has $L^p-L^2$ off-diagonal decay. For $2 \le q < p_+(A)$, $(\Gamma_A(t,s))$ has $L^2-L^q$ off-diagonal decay.
    \end{enumerate}
\end{prop}

\begin{proof}
    First consider \ref{item:p-+(A)-bd}. Since $(\Gamma_A(t,s))$ is a family of contractions on $L^2$, particularly, it is uniformly bounded on $L^2$. By interpolation, it is clear that all the $p \in [1,\infty]$ for which $(\Gamma_A(t,s))$ is uniformly bounded on $L^p$ form an interval containing 2. Let $p_-(A)$ and $p_+(A)$ be the left and right extremes of this interval. It has been shown in \cite[Theorem 1.6]{Zaton2020wp} that there exists $\epsilon>0$ only depending on $n$ and the ellipticity of $A$ so that $(\Gamma_A(t,s))$ is uniformly bounded on $L^p$ for $2-\epsilon<p<2+\epsilon$. We thus infer that $p_-(A)<2$ and $p_+(A)>2$. This proves \ref{item:p-+(A)-bd}.

    The second statement \ref{item:p-+(A)-odd} combines \cite[Lemmas 4.9 and 4.11]{Auscher-Monniaux-Portal2019Lp}. This completes the proof.
\end{proof}

For time-independent coefficients, the critical numbers $p_\pm(A)$ coincide with $p_\pm(L)$ introduced in \cite[\S 3.2]{Auscher2007memoire}, where $L:=-\Div(A(x)\nabla)$. We also know the inequalities $p_-(L)<\frac{2n}{n+2}$ and $p_+(L)>\frac{2n}{n-2}$ are best possible. However, for time-dependent coefficients, we do not know whether the bounds $p_-(A)<2$ and $p_+(A)>2$ are sharp.
\section{Inhomogeneous Cauchy problem}
\label{sec:ic}

This section is concerned with the existence of weak solutions to the inhomogeneous Cauchy problem
\begin{equation*}
    \tag{IC}
    \begin{cases}
    \partial_t v - \Div(A\nabla v) = f \\
    v(0)=0
    \end{cases}
\end{equation*}
as announced in Theorem \ref{thm:wp-ic}. The weak solutions are constructed by extension of the Duhamel operator $\cL^A_1$ defined in \eqref{e:LA1}. We collect the properties of the extension in the following theorem as a general result. Part of the proof is deferred to the end of next section.

\begin{theorem}[Extension of $\cL^A_1$]
    \label{thm:ic}
    Let $\beta>-1/2$ and $p_A(\beta)<p \le \infty$. Then $\cL^A_1$ extends to a bounded operator from $T^p_\beta$ to $T^p_{\beta+1}$, also denoted by $\cL^A_1$. Moreover, the following properties hold for any $f \in T^p_\beta$ and $v:=\cL^A_1(f)$.
    \begin{enumerate}[label=\normalfont(\alph*)]
        \item \label{item:ic_reg} {\normalfont (Regularity)} 
        $v$ lies in $T^{p}_{\beta+1}$ and $\nabla v$ lies in $T^{p}_{\beta+1/2}$ with
        \[ \|v\|_{T^{p}_{\beta+1}} \lesssim \|f\|_{T^{p}_\beta}, \quad \|\nabla v\|_{T^{p}_{\beta+1/2}} \lesssim \|f\|_{T^{p}_\beta}. \]
        
        \item \label{item:ic_sol}
        $v$ is a global weak solution to $\partial_t v - \Div(A\nabla v) = f$.

        \item \label{item:ic-formula} {\normalfont (Explicit formula)}
        It holds that
        \begin{equation}
            \label{e:LA1=LI1+RI1/2}
            v = \cL^\bI_1(f) + \cR^{\bI}_{1/2}((A-\bI)\nabla v).
        \end{equation}
        Here, $\cR^{\bI}_{1/2}$ denotes the extension of the Lions operator to $T^p_{\beta+1/2}$ given by Theorem \ref{thm:lions}.

        \item \label{item:ic_cont} {\normalfont (Continuity and trace)} 
        $v \in C([0,\infty);\scrS')$ with $v(0)=0$. As $t \to 0$, the convergence also occurs in
        \begin{equation}
            \label{e:ic-trace}
            \begin{cases} 
                L^p  & \text{ if } p_A(\beta) < p \le 2 \\
                E^q_\delta \quad & \text{ if } 2<p \le \infty
        \end{cases},
        \end{equation}
        where $\delta>0$ and $q \in [p,\infty]$ are arbitrary parameters.
    \end{enumerate}
    Consequently, $v$ is a global weak solution to \eqref{e:ic} with source term $f$.
\end{theorem}

\begin{remark}
    \label{rem:ic-real-max-reg}
    For time-independent coefficients, the range of $\beta$ and $p$ in Theorem \ref{thm:ic} coincides with that in \cite[Theorem 2]{Auscher-Hou2025_WPMRTent}. They also establish the maximal regularity estimates. Namely, both $\partial_t u$ and $\Div(A\nabla u)$ lie in $T^p_\beta$ with $\|\partial_t u\|_{T^p_\beta} + \|\Div(A\nabla u)\|_{T^p_\beta} \lesssim \|f\|_{T^p_\beta}$. For time-dependent coefficients, such estimates are not clear, since we do not have appropriate estimates for the operator $\varphi \mapsto \Div(A\nabla \Gamma_A(t,s) \varphi)$.
\end{remark}

In this section, we prove \ref{item:ic_reg} and \ref{item:ic_sol}. The proof of \ref{item:ic-formula} and \ref{item:ic_cont} is postponed to Section \ref{ssec:pf-ic-cont}. We employ the theory of singular integral operators on tent spaces developed in \cite{Auscher-Hou2025_WPMRTent}. The following lemma verifies that both $\cL^A_1$ and the backward Duhamel's operator $(\cL^A_1)^\ast$ defined in \eqref{e:LA1*} are involved in this theory. 
\begin{lemma}
    \label{lem:LA-SIO}
    The operator $\cL^A_1$ (resp. $(\cL^A_1)^\ast$) belongs to $\sio^{1+}_{2,q,\infty}$ (resp. $\sio^{1+}_{2,q',\infty}$) for $p_-(A)<q<p_+(A)$.
\end{lemma}

\begin{proof}
    To prove $\cL^A_1 \in \sio^{1+}_{2,q,\infty}$, it suffices to show the kernel $K(t,s):=\I_{\{t>s\}}(t,s)\Gamma_A(t,s)$ belongs to $\sk^1_{2,q,\infty}$, see \cite[Definitions 2 and 3]{Auscher-Hou2025_WPMRTent}. 
    
    To this end, we first show $K$ is strongly measurable, \textit{i.e.}, for any $\phi \in L^2$, the function $(t,s) \mapsto K(t,s)\phi$ is strongly measurable, valued in $L^2$. Let $\psi \in L^2$. We have shown in \eqref{e:s->Gamma(t,s)-weak-meas} that $s \mapsto \langle \Gamma_A(t,s)\phi, \psi \rangle$ is continuous, and Corollary \ref{cor:Gamma} implies $t \mapsto \Gamma_A(t,s)\phi$ is strongly continuous valued in $L^2$, so $t \mapsto \langle \Gamma_A(t,s)\phi, \psi \rangle$ is continuous. Thus, $(t,s) \mapsto \langle \Gamma_A(t,s)\phi, \psi \rangle$ is separately continuous on $\{t>s\}$. This implies $(t,s) \mapsto \Gamma_A(t,s)\phi = K(t,s)\phi$ is weakly measurable, and hence strongly measurable thanks to Pettis' measurability theorem and the fact that $L^2$ is separable. Next, the uniform boundedness of $K(t,s)$ on $L^2$ follows from the fact that $(\Gamma_A(t,s))$ is a family of contractions on $L^2$. Finally, the off-diagonal decay of $K$ comes from Proposition \ref{prop:Gamma-Lp} \ref{item:p-+(A)-odd}. This proves $K \in \sk^1_{2,q,\infty}$ and concludes for $\cL^A_1$. 
    
    For $(\cL^A_1)^\ast$, it follows by duality, see \cite[Proposition 2]{Auscher-Hou2025_WPMRTent}.
\end{proof}

Let us prove Theorem \ref{thm:ic} \ref{item:ic_reg} and \ref{item:ic_sol}.

\begin{proof}[Proof of Theorem \ref{thm:ic} \ref{item:ic_reg} and \ref{item:ic_sol}]
    First consider \ref{item:ic_reg}. The bounded extension of $\cL^A_1$ is a direct consequence of \cite[Proposition 3 and Corollary 2]{Auscher-Hou2025_WPMRTent}, where the conditions are verified in Lemma \ref{lem:LA-SIO}. To prove the gradient estimates, we first pick $f \in L^2_\ssfc(\bR^{1+n}_+)$. Proposition \ref{prop:formula} \ref{item:L2-ic} says $v=\cL^A_1(f)$ is a global weak solution to $\partial_t v - \Div(A\nabla v) = f$. So Corollary \ref{cor:Caccioppoli-tent} yields
    \[ \|\nabla v\|_{T^p_{\beta+1/2}} \lesssim \|v\|_{T^p_{\beta+1}} + \|f\|_{T^p_\beta} \lesssim \|f\|_{T^p_\beta}. \]
    Then a density argument extends the above estimates to all $f \in T^p_{\beta}$ (or weak*-density if $p=\infty$). This proves \ref{item:ic_reg}.

    For \ref{item:ic_sol}, we know from \ref{item:ic_reg} that $v:=\cL^A_1(f)$ lies in $L^2_{\loc}((0,\infty);W^{1,2}_{\loc})$, since all the tent spaces embed into $L^2_{\loc}(\bR^{1+n}_+)$. Moreover, for $f \in L^2_{\ssfc}(\bR^{1+n}_+)$, $v$ satisfies $\partial_t v - \Div(A\nabla v) = f$ in $\scrD'(\bR^{1+n}_+)$. Using the boundedness of $\cL^A_1$ and $\nabla \cL^A_1$, one can extend this identity (valued in $\scrD'(\bR^{1+n}_+)$) to all $f \in T^p_\beta$ by density, or weak*-density if $p=\infty$.
\end{proof}

\begin{cor}
    \label{cor:LA1*}
    Let $\beta>-1/2$ and $1 \le p' < \max\{p_A(\beta),1\}'$. Then $(\cL^A_1)^\ast$ extends to a bounded operator from $T^{p'}_{-\beta-1}$ to $T^{p'}_{-\beta}$, also denoted by $(\cL^A_1)^\ast$. It additionally satisfies
    \[ \|\nabla (\cL^A_1)^\ast(f) \|_{T^{p'}_{-\beta-1/2}} \lesssim \|f\|_{T^{p'}_{-\beta-1}}. \]
\end{cor}

\begin{proof}
    The statements follow from adapting the arguments in the proof of Theorem \ref{thm:ic} \ref{item:ic_reg} to backward singular integral operators (see \cite[Proposition 4 and Corollary 3]{Auscher-Hou2025_WPMRTent}) and backward equations. The conditions are also verified in Lemma \ref{lem:LA-SIO}. We only need to mention Proposition \ref{prop:formula} \ref{item:L2-ic-back} says that for $f \in L^2_\ssfc(\bR^{1+n}_+)$, $w:=(\cL^A_1)^\ast(f)$ is a weak solution to the backward equation $-\partial_s w - \Div(A(s)^\ast \nabla w) = f$. Detailed verification is left to the reader.
\end{proof}
\section{Lions' equation}
\label{sec:lions}

In this section, we construct weak solutions to the \emph{Lions equation}
\begin{equation}
    \label{e:lions}
    \tag{L}
    \begin{cases}
        \partial_t u - \Div(A\nabla u) = \Div F \\
        u(0)=0
    \end{cases}
\end{equation}
as asserted in Theorem \ref{thm:wp-lions}, by extension of the Lions operator $\cR^A_{1/2}$ defined in \eqref{e:RA1/2}. The main theorem of this section summarizes the properties of the extension. We use the critical exponents $p^\pm_\zeta(\beta)$ defined in \eqref{e:p-zeta(beta)} and \eqref{e:p+zeta(beta)}.

\begin{theorem}[Extension of $\cR^A_{1/2}$]
    \label{thm:lions}
    There exists $\beta_A \in [-1,-1/2)$ only depending on the ellipticity of $A$ and the dimension $n$ so that, for $\beta>\beta_A$ and
    \[ \begin{cases}
        p^-_{\beta_A}(\beta) < p \le p^+_{\beta_A}(\beta) = \infty & \text{ if } \beta \ge -1/2 \\
        p^-_{\beta_A}(\beta) < p < p^+_{\beta_A}(\beta) & \text{ if } \beta_A < \beta < -1/2 
    \end{cases}, \]
    the Lions operator $\cR^A_{1/2}$ extends to a bounded operator from $T^p_{\beta+1/2}$ to $T^p_{\beta+1}$, also denoted by $\cR^A_{1/2}$. Moreover, the following properties hold for any $F \in T^p_{\beta+{1}/{2}}$ and $u:=\cR^A_{1/2}(F)$.
    \begin{enumerate}[label=\normalfont(\alph*)]
        \item \label{item:lions_reg} {\normalfont (Regularity)}
        $u$ lies in $T^p_{\beta+1}$ and $\nabla u$ lies in $T^p_{\beta+1/2}$ with 
        \[ \|u\|_{T^{p}_{\beta+1}} \lesssim \|F\|_{T^{p}_{\beta+{1}/{2}}}, \quad \|\nabla u\|_{T^{p}_{\beta+{1}/{2}}} \lesssim \|F\|_{T^{p}_{\beta+{1}/{2}}}. \]

        \item \label{item:lions_sol}
        $u$ is a global weak solution to $\partial_t u - \Div(A\nabla u) = \Div F$.

        \item \label{item:lions_formula} {\normalfont (Explicit formula)}
        Define $\tilF:=(A-\bI)\nabla u+F$. Then 
        \begin{equation}
            \label{e:lions-formula}
            u = \cR^{\bI}_{1/2}(\tilF) = \Div \cL_1^{\bI}(\tilF) \quad \text{ in } \scrD'(\bR^{1+n}_+).
        \end{equation}

        \item \label{item:lions_cont} {\normalfont (Continuity and traces)}
        $u \in C([0,\infty);\scrS')$ with $u(0)=0$. As $t \to 0$, the convergence also occurs in the following trace spaces shown in Table \ref{tab:lions_cont}, with arbitrary parameters $\delta>0$, $q \in [p,\infty]$, and $s \in [-1,2\beta+1]$.

        \begin{table}[!ht]
            \centering
            \caption{Trace spaces for $\cR^A_{1/2}(F)$.}
            \label{tab:lions_cont}
            \begin{tabular}{c|c|c|c}
            Conditions 
            & $p^-_{\beta_A}(\beta) <p \le 2$ 
            & $2<p<\infty$
            & $p=\infty$ \\
            \hline
         
            $\beta>-1/2$ 
            & $L^p$ 
            & $E^q_\delta \cap E^{-1,q}_\delta$
            & $E^\infty_\delta \cap E^{-1,\infty}_\delta$ \\
            \hline
        
            $\beta = -1/2$ 
            & $L^p$
            & $E^q_\delta \cap E^{-1,q}_\delta$
            & $L^2_{\loc} \cap E^{-1,\infty}_\delta$ \\
            \hline

            $\beta_A<\beta<-1/2$ 
            & $\DotH^{s,p}$
            & $E^{-1,q}_\delta$
            & not applicable
            \end{tabular}
        \end{table}
    \end{enumerate}
    Consequently, $u$ is a global weak solution to the Lions equation \eqref{e:lions}.
\end{theorem}

\begin{remark}
    \label{rem:beta<-1/2-A(x)}
    For time-independent coefficients, it has been shown in \cite[Theorem 5.1]{Auscher-Hou2024-HCL} that $\beta_A=-1$, and the range of $p$ for which these properties are valid strictly contains $p^-_{-1}(\beta)<p \le \infty$. But here, we present new trace spaces $E^q_\delta$ for $\beta \ge -1/2$ and $2<p \le \infty$, which will be useful for proving uniqueness.
\end{remark}

\begin{remark}
    \label{rem:lions-endpt-Einfty-bd}
    As we shall see in the proof, for $\beta=-1/2$ and $p=\infty$, we still obtain that for any $\delta>0$, $\sup_{0<t<\delta} \|u(t)\|_{E^\infty_\delta}$ is uniformly bounded, although we do not know how to prove the trace in $E^\infty_\delta$.
\end{remark}

\begin{remark}
    \label{rem:lions-RA1/2->abcde'}
    In the range of $\beta>-1$ and $\frac{n}{n+2\beta+2} < p \le \infty$, once one can show that $\cR^A_{1/2}$ extends to a bounded operator from $T^p_{\beta+1/2}$ to $T^p_{\beta+1}$, then all subsequent properties can be proved by the same arguments. Moreover, for $p=\infty$, if $\cR^A_{1/2}$ extends to a bounded operator from $T^\infty_{\beta+1/2}$ to $T^\infty_{\beta+1}$ for some $\beta \in (-1,-1/2)$, then for any $F \in T^\infty_{\beta+1/2}$, $\cR^A_{1/2}(F)(t)$ tends to 0 in $E^{-1,\infty}_{\delta}$ as $t \to 0$.
\end{remark}

Let us mention a very recent work \cite{Auscher-Portal2025-Lions}, which addresses the case $\beta=-1/2$ and $2 \le p \le \infty$. For real coefficients, their results extend to $1-\epsilon<p<2$ by exploiting local H\"older continuity of weak solutions. However, this continuity property fails for complex coefficients, even in the autonomous case. Instead, our approach furnishes the case $p_-(A)<p<2$ via an extrapolation argument and further includes new cases $\beta \ne -1/2$. We also establish the time continuity of weak solutions.

The proof of Theorem \ref{thm:lions} is presented in Section \ref{ssec:pf-thm-lions}. We first prove three lemmas on the bounded extension of $\cR^A_{1/2}$ in different ranges. 
\begin{lemma}
    \label{lem:lions-beta>-1/2}
    Let $\beta>-1/2$ and $p_A(\beta)<p \le \infty$. Then $\cR^A_{1/2}$ extends to a bounded operator from $T^p_{\beta+1/2}$ to $T^p_{\beta+1}$.
\end{lemma}

\begin{lemma}
    \label{lem:lions-beta<-1/2-p=2}
    There exists $\epsilon>0$ only depending on the ellipticity of $A$ and the dimension $n$ so that the properties of $\cR^A_{1/2}$ in Lemma \ref{lem:lions-beta>-1/2} are also valid for $-1/2-\epsilon<\beta<-1/2+\epsilon$ and $p=2$.
\end{lemma}

\begin{lemma}
    \label{lem:lions-beta=-1/2-p=infty}
    The properties of $\cR^A_{1/2}$ in Lemma \ref{lem:lions-beta>-1/2} are also valid for $\beta=-1/2$ and $p=\infty$.
\end{lemma}

The proofs are provided right below in Sections \ref{ssec:pf-lions-beta>-1/2}, \ref{ssec:pf-lions-beta<-1/2-p=2}, and \ref{ssec:pf-thm-lions}. 

\subsection{Proof of Lemma \ref{lem:lions-beta>-1/2}}
\label{ssec:pf-lions-beta>-1/2}
Let $\beta>-1/2$ and $p_A(\beta)<p \le \infty$. We split the proof into two cases. \\

\paragraph{Case 1: $\max\{p_A(\beta),1\} < p \le \infty$}
We argue by duality. Let $F \in L^2_\ssfc(\bR^{1+n}_+)$ and $\psi \in \Cc(\bR^{1+n}_+)$. Fubini's theorem yields
\begin{align*}
    \langle \cR^A_{1/2}(F),\psi \rangle_{L^2(\bR^{1+n}_+)} 
    &= - \int_0^\infty \int_{\bR^n} F(s,y) \left( \int_s^{\infty} \overline{\nabla \Gamma_A(t,s)^\ast \psi(t)}(y) dt \right) dsdy \\
    &= - \langle F,\nabla (\cL^A_1)^\ast(\psi) \rangle_{L^2(\bR^{1+n}_+)}.
\end{align*} 
We further apply duality of tent spaces and Corollary \ref{cor:LA1*} to get
\[ | \langle \cR^A_{1/2}(F),\psi \rangle_{L^2(\bR^{1+n}_+)} | \lesssim \|F\|_{T^p_{\beta+1/2}} \|\nabla (\cL^A_1)^\ast \psi\|_{T^{p'}_{-\beta-1/2}} \lesssim \|F\|_{T^p_{\beta+1/2}} \|\psi\|_{T^{p'}_{-\beta-1}}. \]
The bounded extension of $\cR^A_{1/2}$ hence follows by density, or weak*-density if $p=\infty$. \\

\paragraph{Case 2: $p_A(\beta)<p \le 1$}
We use atomic decomposition of tent spaces, see \cite[Proposition 5]{Coifman-Meyer-Stein1985_TentSpaces}. Recall that for $\beta \in \bR$ and $0<p \le 1$, a function $a \in L^2_{\beta+1/2}(\bR^{1+n}_+)$ is called a \emph{$T^p_{\beta+1/2}$-atom}, if there exists a ball $B \subset \bR^n$ so that $\supp(a) \subset [0,r(B)^2] \times B$, and
\[ \|a\|_{L^2_{\beta+1/2}(\bR^{1+n}_+)} \le |B|^{-[p,2]}. \]
The ball $B$ is called \emph{associated} to $a$. Write $r:=r(B)$, $C_0:=2B$, and $C_j:=2^{j+1} B \setminus 2^j B$ for $j \ge 1$. For $j \ge 4$, define
\begin{align*}
    M^{(1)}_j &:= \left( 0,(2^3 r)^2 \right] \times C_j, \\
    M^{(2)}_j &:= \left( (2^3 r)^2,(2^j r)^2 \right) \times C_j, \\
    M^{(3)}_j &:= \left[ (2^j r)^2,(2^{j+1} r)^2 \right) \times 2^{j+1} B.
\end{align*}

The following lemma describes the molecular decay of $\cR^A_{1/2}$ acting on atoms.
\begin{lemma}[Molecular decay]
    \label{lem:lions-mol}
    Let $\beta>-1/2$, $0<p \le 1$, and $p_-(A)<q<2$. There exists a constant $c>0$ depending on $A$ and $q$ so that for any $T^p_{\beta+1/2}$-atom $a$ with an associated ball $B \subset \bR^n$, the following estimates hold for $u:=\cR^A_{1/2}(a)$ and $j \ge 4$,
    \begin{align}
        &\label{e:lions_atom_Mj1}
        \| \I_{M_j^{(1)}} u\|_{ L^2_{\beta+1}(\bR^{1+n}_+) } \lesssim 2^{jn[p,2]} e^{-c2^{2j}} |2^{j+1} B|^{-[p,2]},  \\
        &\label{e:lions_atom_Mj2}
        \| \I_{M_j^{(2)}} u\|_{ L^2_{\beta+1}(\bR^{1+n}_+) }  \lesssim 2^{-j(2\beta+1+n[q,p])} |2^{j+1} B|^{-[p,2]}, \\
        &\label{e:lions_atom_Mj3}
        \| \I_{M_j^{(3)}} u\|_{ L^2_{\beta+1}(\bR^{1+n}_+) } \lesssim 2^{-j(2\beta+1+n[q,p])} |2^{j+1} B|^{-[p,2]}.
    \end{align}
\end{lemma}

The proof is provided right below. Admitting this lemma, let us prove the bounded extension of $\cR^A_{1/2}$. Thanks to \cite[Lemma 6]{Auscher-Hou2025_WPMRTent}, it suffices to verify $\cR^A_{1/2}$ is uniformly bounded on $T^p_{\beta+1/2}$-atoms. Let $a$ be a $T^p_{\beta+1/2}$-atom, $B \subset \bR^n$ be a ball associated to $a$, and $u:=\cR^A_{1/2}(a)$. Write $Q_0:= \left( 0,(2^3 r(B))^2 \right) \times 2^3 B$. We get from Case 1 that
\begin{equation}
    \label{e:lions-1Q0-u}
    \|\I_{Q_0} u\|_{L^2_{\beta+1}(\bR^{1+n}_+)} \lesssim \|a\|_{L^2_{\beta+1/2}(\bR^{1+n}_+)} \lesssim |2^3 B|^{-[p,2]}.
\end{equation}
Meanwhile, pick $q<2$ sufficiently close to $p_-(A)$ so that
\[ p_A(\beta) = \frac{np_-(A)}{n+(2\beta+1)p_-(A)} < \frac{nq}{n+(2\beta+1)q} < p \le 1. \]
In particular, we get $2\beta+1+n[q,p]>0$. Applying Lemma \ref{lem:lions-mol} to such $q$ gives estimates of $u$ on molecules $M^{(i)}_j$ for $1 \le i \le 3$ and $j \ge 4$. These estimates, together with \eqref{e:lions-1Q0-u}, yield
\[ \|u\|_{T^p_{\beta+1}}^p \lesssim 1+\sum_{j \ge 4} 2^{jn[p,2]p} e^{-c2^{2j}p} + \sum_{j \ge 4} 2^{-j(2\beta+1+n[q,p])p} \le C < \infty. \] 
This completes the proof of Lemma \ref{lem:lions-beta>-1/2}.

\begin{proof}[Proof of Lemma \ref{lem:lions-mol}]
    The proof of \eqref{e:lions_atom_Mj2} and \eqref{e:lions_atom_Mj3} is a verbatim adaptation of \cite[Lemma 5.9, (5.4) and (5.5)]{Auscher-Hou2024-HCL}, using the $L^2-L^{q'}$ off-diagonal decay of $(\Gamma_A(t,s)^\ast)$. Let us focus on \eqref{e:lions_atom_Mj1}. As $a \in L^2_{\beta+1/2}(\bR^{1+n}_+)$, we show in Case 1 that $u:=\cR^A_{1/2}(a)$ lies in $L^2_{\beta+1}(\bR^{1+n}_+)$. Let $\phi \in \Cc(\bR^n)$ with $\supp(\phi) \subset C_j$. We have
    \[ \langle u(t),\phi \rangle_{L^2(\bR^n)} = -\int_0^t \int_{\bR^n} a(s,y) \overline{\nabla \Gamma_A(t,s)^\ast \phi}(y) dsdy. \]
    Write $v_t(s,y):=\Gamma_A(t,s)^\ast \phi(y)$. As $2\beta+1>0$, applying the properties of $a$ and Cauchy--Schwarz inequality gives
    \begin{align*}
        \left| \langle u(t),\phi \rangle_{L^2(\bR^n)} \right|^2 
        &\le \|a\|_{L^2_{\beta+1/2}(\bR^{1+n}_+)}^2 \int_0^t \int_B s^{2\beta+1} |\nabla v_t(s,y)|^2 dsdy \\
        &\le |B|^{-2[p,2]} t^{2\beta+1} \int_0^t \int_B |\nabla v_t(s,y)|^2 dsdy.
    \end{align*}
    Then we take a covering of $(0,t) \times B$ by reverse Whitney cubes $( (1-2^{-k})t,(1-2^{-k-1})t ) \times B( x_0,(1-2^{-k-1})^{1/2}t^{1/2} )$ for $k \ge 0$. Corollary \ref{cor:Gamma-} implies $v_t$ is a weak solution to the backward equation $-\partial_s v_t - \Div (A(s)^\ast \nabla v_t) = 0$. By applying (backward) Caccioppoli's inequality to each Whitney cube and summing the results, we obtain
    \[ \int_0^t \int_B |\nabla v_t(s,y)|^2 dsdy \lesssim \int_0^t \int_{15B} \frac{|v_t(s,y)|^2}{t-s} dsdy. \]
    The $L^2-L^2$ off-diagonal estimates of $(\Gamma_A(t,s)^\ast)$ provide us with a constant $c>0$ so that
    \[ \int_0^t \int_{15B} \frac{|v_t(s,y)|^2}{t-s} dsdy \lesssim \int_0^t \frac{1}{t-s} e^{-\frac{2c(2^j r)^2}{t-s}} \|\phi\|_{L^2}^2 ds \lesssim e^{-\frac{c (2^j r)^2}{t}} \|\phi\|_{L^2}^2. \]
    Gathering these estimates, by duality, we obtain
    \[ \|u(t)\|_{L^2(C_j)}^2 \lesssim |B|^{-2[p,2]} t^{2\beta+1} e^{-\frac{c (2^j r)^2}{t}}, \quad 0<t \le (2^3 r)^2. \]
    Integrating both sides over $0<t \le (2^3 r)^2$ gives \eqref{e:lions_atom_Mj1} as desired. 
\end{proof}

\subsection{Proof of Lemma \ref{lem:lions-beta<-1/2-p=2}}
\label{ssec:pf-lions-beta<-1/2-p=2}
We use an abstract argument relying on Sneiberg's lemma. Let $\beta \in \bR$ and $1<p<\infty$. Define
\[ S^p_{\beta+1/2} := \left\{ u \in L^2_{\loc}((0,\infty);W^{1,2}_{\loc}) : \nabla u \in T^p_{\beta+1/2}, ~ \partial_t u \in \Div T^p_{\beta+1/2} \right\}, \]
endowed with the semi-norm
\[ \|u\|_{S^p_{\beta+1/2}} := \|\nabla u\|_{T^p_{\beta+1/2}} + \|\partial_t u\|_{\Div T^p_{\beta+1/2}}. \]

\begin{lemma}
    \label{lem:tr-Sp}
    Let $-1<\beta<0$ and $1<p<\infty$. Let $u \in S^p_{\beta+1/2}$. Then $u$ belongs to $C([0,\infty);\scrS')$. Moreover, $u(0)$ belongs to $\DotH^{2\beta+1,p}+\bC$ with
    \[ \|u(0)\|_{\DotH^{2\beta+1,p}/\bC} \lesssim \|u\|_{S^p_{\beta+1/2}}. \]
    We call $u(0)$ the \emph{trace} of $u$, denoted by $\tr(u)$.
\end{lemma}

\begin{proof}
    Let $u \in S^p_{\beta+1/2}$. As $\partial_t u - \Delta u \in \Div T^p_{\beta+1/2}$, \cite[Theorem 1.6]{Auscher-Hou2024-HCL} says there is a unique weak solution $v \in S^p_{\beta+1/2}$ to the Cauchy problem
    \[ \begin{cases}
        \partial_t v - \Delta v = \partial_t u - \Delta u \\
        v(0)=0
    \end{cases}. \]
    It additionally belongs to $C([0,\infty);\scrS')$ and satisfies
    \[ \|v\|_{S^p_{\beta+1/2}} \lesssim \|\partial_t u - \Delta u\|_{\Div T^p_{\beta+1/2}} \lesssim \|u\|_{S^p_{\beta+1/2}}. \]
    Observe that $w:=u-v$ lies in $S^p_{\beta+1/2}$ and $w$ is a weak solution to the heat equation $\partial_t w - \Delta w = 0$. Then we invoke \cite[Theorem 1.1]{Auscher-Hou2024-HCL} to get that there is a unique $u_0 \in \scrS'$ so that $w=\cE_\bI(u_0)$, and
    \[ \|u_0\|_{\DotH^{2\beta+1,p}/\bC} \eqsim \|\nabla w\|_{T^p_{\beta+1/2}} \le \|\nabla u\|_{T^p_{\beta+1/2}} + \|\nabla v\|_{T^p_{\beta+1/2}} \lesssim \|u\|_{S^p_{\beta+1/2}}. \]
    Note that for any $t>0$, $w(t)=\cE_\bI(u_0)(t) = e^{t\Delta} u_0$, so $w$ also belongs to $C([0,\infty);\scrS')$ with $w(0)=u_0$. Therefore, we get $u \in C([0,\infty);\scrS')$ with $u(0)=u_0$. This completes the proof.
\end{proof}

\begin{prop}
    \label{prop:Sp-iso}
    Let $-1<\beta<0$ and $1<p<\infty$. The map
    \begin{align*}
        \Phi : S^p_{\beta+1/2} &\to \Div T^p_{\beta+1/2} \times (\DotH^{2\beta+1,p}+\bC) \\
        \Phi(u) &:= (\partial_t u -\Delta u, \qquad \tr(u)).
    \end{align*}
    is an isomorphism of semi-normed spaces. Hence, the quotient map of $\Phi$ (induced by the canonical projection modulo constants) is an isomorphism of Banach spaces from $S^p_{\beta+1/2}/\bC$ to $\Div T^p_{\beta+1/2} \times \DotH^{2\beta+1,p}/\bC$.
\end{prop}

\begin{proof}
    Lemma \ref{lem:tr-Sp} shows $\Phi$ is bounded. Let us construct the inverse $\Psi$. For any $g \in \Div T^p_{\beta+1/2}$ and $v_0 \in \DotH^{2\beta+1,p} + \bC$, let $v=:\Psi(g,v_0)$ be the unique weak solution in $S^p_{\beta+1/2}$ to the Cauchy problem
    \begin{equation}
        \label{e:Psi(g,v0)-eq}
        \begin{cases}
            \partial_t v - \Delta v = g \\
            v(0)=v_0
        \end{cases}.
    \end{equation}
    The existence and uniqueness are ensured by \cite[Theorem 1.7]{Auscher-Hou2024-HCL}. We also infer that
    \[ \|\Psi(g,v_0)\|_{S^p_{\beta+1/2}} \lesssim \|g\|_{\Div T^p_{\beta+1/2}} + \|v_0\|_{\DotH^{2\beta+1,p}/\bC}, \]
    so $\Psi$ is bounded. Let us verify that $\Psi$ is the inverse of $\Phi$. First, for any $g \in \Div T^p_{\beta+1/2}$ and $v_0 \in \DotH^{2\beta+1,p}+\bC$, the equation \eqref{e:Psi(g,v0)-eq} shows $\partial_t \Psi(g,v_0) - \Delta \Psi(g,v_0) = g$, and $\Psi(g,v_0)(t)$ converges to $v_0$ as $t \to 0$ in $\scrD'$, so $v_0$ coincides with $\tr(\Psi(g,v_0))$ constructed in Lemma \ref{lem:tr-Sp}. We hence get $\Phi \circ \Psi(g,v_0) = (g,v_0)$.
    
    On the other hand, for any $u \in S^p_{\beta+1/2}$, note that both $u$ and $w:=\Psi \circ \Phi(u)$ are weak solutions in $S^p_{\beta+1/2}$ to the Cauchy problem
    \begin{equation}
        \label{e:Psi-Phi-u-eq}
        \begin{cases}
            \partial_t w - \Delta w = \partial_t u - \Delta u \\
            w(0)=\tr(u)
        \end{cases}.
    \end{equation}
    Thanks to the uniqueness of weak solutions to \eqref{e:Psi-Phi-u-eq} in $S^p_{\beta+1/2}$, we get $u = w = \Psi \circ \Phi(u)$. This proves $\Psi$ is the inverse of $\Phi$.

    The last point follows from the fact that $\Phi(c)=(0,c)$ for any constant function $c$. This completes the proof.
\end{proof}

Let us finish the proof of Lemma \ref{lem:lions-beta<-1/2-p=2}.
\begin{proof}[Proof of Lemma \ref{lem:lions-beta<-1/2-p=2}]
    Let $-1<\beta<0$ and $p=2$. Define the map
    \begin{equation}
        \label{e:PhiA}
        \begin{aligned}
            \Phi_A : S^2_{\beta+1/2} &\to \Div T^2_{\beta+1/2} \times (\DotH^{2\beta+1,2}+\bC) \\
            \Phi_A(u) &:= (\partial_t u -\Div(A\nabla u), \quad \tr(u)).
        \end{aligned}
    \end{equation}
    
    We first show $\Phi_A$ is an isomorphism of semi-normed spaces for $\beta=-1/2$. Indeed, let $g \in \Div T^2_0$ and $v_0 \in \DotH^{0,2}+\bC \simeq L^2 + \bC$. Let $v=:\Psi_A(g,v_0)$ be the unique weak solution in $S^2_0$ to the Cauchy problem
    \begin{equation}
        \label{e:PsiA(g,v0)}
        \begin{cases}
            \partial_t v - \Div(A\nabla v) = g \\
            v(0)=v_0
        \end{cases}
    \end{equation}
    given by Proposition \ref{prop:L2}. It also satisfies $\|v\|_{S^2_0} \lesssim \|g\|_{\Div T^2_0} + \|v_0\|_{\DotH^{0,2}/\bC}$, so $\Psi_A$ is bounded. Ensured by the uniqueness of weak solutions to \eqref{e:PsiA(g,v0)} in $S^2_0$, one can get $\Psi_A$ is the inverse of $\Phi_A$ by adapting the proof of Proposition \ref{prop:Sp-iso} \textit{mutatis mutandis}.

    Next, observe that for any $F \in T^2_0$,
    \begin{equation}
        \label{e:PhiA-1=RA1/2}
        \Phi_A^{-1}(\Div F,0) = \cR^A_{1/2}(F).
    \end{equation}
    Indeed, both $u:=\Phi_A^{-1}(\Div F,0)$ and $v:=\cR^A_{1/2}(F)$ are weak solutions in $S^2_0$ to the Cauchy problem
    \begin{equation}
        \label{e:PhiA-1(divF,0)-eq}
        \begin{cases}
            \partial_t u - \Div (A\nabla u) = \Div F \\
            u(0)=0
        \end{cases}.
    \end{equation}
    Then \eqref{e:PhiA-1=RA1/2} follows from uniqueness of weak solutions to \eqref{e:PhiA-1(divF,0)-eq} in $S^2_0$.

    We claim that there is a constant $\epsilon>0$ only depending on the ellipticity of $A$ and the dimension $n$ so that $\Phi_A$ is an isomorphism of semi-normed spaces for $-1/2-\epsilon<\beta<-1/2+\epsilon$. Suppose it holds. Then the map $F \mapsto \nabla \Phi_A^{-1}(\Div F,0)$ extends to $T^2_{\beta+1/2}$ with
    \[ \|\nabla \Phi_A^{-1}(\Div F,0)\|_{T^2_{\beta+1/2}} \le \|\Phi_A^{-1}(\Div F,0)\|_{S^2_{\beta+1/2}} \lesssim \|F\|_{T^2_{\beta+1/2}}. \]
    Moreover, for any $F \in T^2_{\beta+1/2} \cap T^2_0$, applying Proposition \ref{prop:formula} \ref{item:formula-L2-RA1/2} and \eqref{e:PhiA-1=RA1/2} yields
    \[ \cR^A_{1/2}(F) = \cR^{\bI}_{1/2}\left( (A-\bI) \nabla \Phi_A^{-1}(\Div F,0) +F \right). \]
    Using boundedness of $\cR^{\bI}_{1/2}$ from $T^2_{\beta+1/2}$ to $T^2_{\beta+1}$ for $\beta>-1$ (see \cite[Theorem 5.1]{Auscher-Hou2024-HCL}), we obtain
    \[ \|\cR^A_{1/2}(F)\|_{T^2_{\beta+1}} \lesssim \|\nabla \Phi_A^{-1}(\Div F,0)\|_{T^2_{\beta+1/2}} + \|F\|_{T^2_{\beta+1/2}} \lesssim \|F\|_{T^2_{\beta+1/2}}. \] 
    Thus, by density, $\cR^A_{1/2}$ extends to a bounded operator from $T^2_{\beta+1/2}$ to $T^2_{\beta+1}$ for $-1/2-\epsilon<\beta<-1/2+\epsilon$.
    
    It only remains to prove the claim. To this end, we consider the bounded map $\widetilde{\Phi}_A$ induced by $\Phi_A$ on the quotient space via the canonical projection modulo constants, noting that $\Phi_A(c)=(0,c)$ for any constant function $c$. We have the commutative diagram as shown in Figure \ref{fig:Phi_A-diagram}.
    \begin{figure}[ht]
        \centering
        \begin{tikzcd}
            S^2_{\beta+1/2} \arrow[d] \arrow[r, "\Phi_A"] & {\Div T^2_{\beta+1/2} \times (\DotH^{2\beta+1,2}+\bC)} \arrow[d] \\
            S^2_{\beta+1/2}/\bC \arrow[r, "\widetilde{\Phi}_A"]     & {\Div T^2_{\beta+1/2} \times \DotH^{2\beta+1,2}/\bC} 
        \end{tikzcd}
        \caption{Quotient map of $\Phi_A$}
        \label{fig:Phi_A-diagram}
    \end{figure} 
    
    Proposition \ref{prop:Sp-iso} implies $S^2_{\beta+1/2}/\bC$ and $\Div T^2_{\beta+1/2} \times \DotH^{2\beta+1,2}/\bC$ are two isomorphic complex interpolation scales of Banach spaces. We also infer from the first point that $\widetilde{\Phi}_A$ is an isomorphism of Banach spaces for $\beta=-1/2$. Therefore, Sneiberg's Lemma (see \cite{Sneiberg1974-SneibergLemma} or \cite[Theorem 2.7]{Kalton-Mitrea1998-interpolation}) yields there is a constant $\epsilon>0$  so that for $-1/2-\epsilon<\beta<-1/2+\epsilon$, $\widetilde{\Phi}_A$ is an isomorphism of Banach spaces. By lifting, $\Phi_A$ is an isomorphism of semi-normed spaces. This completes the proof.
\end{proof}

\subsection{Proof of Theorem \ref{thm:lions}}
\label{ssec:pf-thm-lions}

Let us first prove Lemma \ref{lem:lions-beta=-1/2-p=infty}.

\begin{proof}[Proof of Lemma \ref{lem:lions-beta=-1/2-p=infty}]
    Define the operator $\cR^A_0$ on $L^2(\bR^{1+n}_+)$ by
    \begin{equation}
        \label{e:RA0}
        \cR^A_0(F) := \nabla \cR^A_{1/2}(F), \quad F \in L^2(\bR^{1+n}_+).
    \end{equation}
    It has been shown in \cite[Theorem 3.1]{Auscher-Portal2025-Lions} that $\cR^A_0$ extends to a bounded operator on $T^\infty_0$. Let $F \in L^2_\ssfc(\bR^{1+n}_+)$ and $\tilF:=(A-\bI)\cR^A_0(F) + F$. Proposition \ref{prop:formula} \ref{item:formula-L2-RA1/2} says $\cR^A_{1/2}(F) = \cR^\bI_{1/2}(\tilF)$. So we get
    \begin{align*}
        \|\cR^A_{1/2}(F)\|_{T^\infty_{1/2}} 
        &= \|\cR^\bI_{1/2}(\tilF)\|_{T^\infty_{1/2}} \lesssim \|\tilF\|_{T^\infty_0} \\
        &\lesssim \|\cR^A_0(F)\|_{T^\infty_0} + \|F\|_{T^\infty_0} \lesssim \|F\|_{T^\infty_0},
    \end{align*}
    where the first inequality comes from boundedness of $\cR^\bI_{1/2}$ from $T^\infty_0$ to $T^\infty_{1/2}$, see \cite[Theorem 5.1]{Auscher-Hou2024-HCL}. It hence implies $\cR^A_{1/2}$ extends to a bounded operator from $T^\infty_0$ to $T^\infty_{1/2}$ by weak*-density.
\end{proof}

Now, we present the proof of Theorem \ref{thm:lions}.

\begin{proof}[Proof of Theorem \ref{thm:lions}]
Let $\epsilon>0$ be the constant given by Lemma \ref{lem:lions-beta<-1/2-p=2}. Define
\begin{equation}
    \label{e:betaA-epsilon}
    \boxed{\beta_A:=-1/2-\epsilon}.
\end{equation}
Let $\beta>\beta_A$ and
\[ \begin{cases}
    p^-_{\beta_A}(\beta) < p \le p^+_{\beta_A}(\beta) = \infty & \text{ if } \beta \ge -1/2 \\
    p^-_{\beta_A}(\beta) < p < p^+_{\beta_A}(\beta) & \text{ if } \beta_A < \beta < -1/2 
\end{cases}. \]

The extension of $\cR^A_{1/2}$ follows by interpolation of tent spaces, thanks to Lemmas \ref{lem:lions-beta>-1/2}, \ref{lem:lions-beta<-1/2-p=2}, and \ref{lem:lions-beta=-1/2-p=infty}. The numbers $p^-_{\beta_A}(\beta)$ and $p^+_{\beta_A}(\beta)$ are exactly designed for this interpolation argument.

Let us verify the properties. First consider \ref{item:lions_reg}. We argue by density. For $F \in L^2_\ssfc(\bR^{1+n}_+)$,  Proposition \ref{prop:formula} \ref{item:L2-HCL} says $u:=\cR^A_{1/2}(F)$ is a global weak solution to $\partial_t u - \Div (A \nabla u) = \Div F$. Then Corollary \ref{cor:Caccioppoli-tent} gives
\[ \|\nabla u\|_{T^p_{\beta+1/2}} \lesssim \|u\|_{T^p_{\beta+1}} + \|F\|_{T^p_{\beta+1/2}} \lesssim \|F\|_{T^p_{\beta+1/2}}. \]
So the gradient estimate holds for $F \in L^2_\ssfc(\bR^{1+n}_+)$. One can extend it to all $F \in T^p_{\beta+1/2}$ by density or weak*-density if $p=\infty$ (we shall omit to mention this subtlety in the sequel). 

For \ref{item:lions_sol}, as $u \in T^p_{\beta+1}$ and $\nabla u \in T^p_{\beta+1/2}$, we get $u \in L^2_{\loc}((0,\infty);W^{1,2}_{\loc})$. For any $F \in L^2_{\ssfc}(\bR^{1+n}_+)$, $u$ satisfies $\partial_t u - \Div(A\nabla u) = \Div F$ in $\scrD'(\bR^{1+n}_+)$. Using the \textit{a priori} estimates in \ref{item:lions_reg}, by density, one can extend the identity to all $F \in T^p_{\beta+1/2}$, valued in $\scrD'(\bR^{1+n}_+)$. This proves \ref{item:lions_sol}.

To prove \ref{item:lions_formula}, we also argue by density. Proposition \ref{prop:formula} \ref{item:formula-L2-RA1/2} shows the formula \eqref{e:lions-formula} holds for $F \in T^p_{\beta+1/2} \cap L^2(\bR^{1+n}_+)$. Note that all the operators involved have bounded extensions to $T^p_{\beta+1/2}$: For $\cR^A_{1/2}$, $\nabla \cR^A_{1/2}$, and $\cR^\bI_{1/2}$, it follows from \ref{item:lions_reg}; For $\cL^{\bI}_1$, it follows from Theorem \ref{thm:ic} \ref{item:ic_reg}, as $\beta+1/2>-1/2$ and $p^-_{\beta_A}(\beta)>\frac{n}{n+2\beta+2}=p_\bI(\beta+1/2)$. Then by density, it holds for all $F \in T^p_{\beta+1/2}$, valued in $L^2_{\loc}(\bR^{1+n}_+)$.

Finally consider \ref{item:lions_cont}. Write $\tilF:=(A-\bI)\nabla u + F$, so $u=\cR^\bI_{1/2}(\tilF)$ (by \ref{item:lions_formula}). As $\beta>\beta_A \ge -1$ and $p>p^-_{\beta_A}(\beta)>\frac{n}{n+2\beta+2}$, \cite[Theorem 5.1(d)]{Auscher-Hou2024-HCL} yields $\cR^\bI_{1/2}(\tilF)$ lies in $C([0,\infty);\scrS')$, and so does $u$. 

The proof of trace spaces is split into three cases. \\

\paragraph{Case 1: $p^-_{\beta_A}(\beta)<p \le 2$}
It follows from the same arguments right above, also using the formula in \ref{item:lions_formula} and \cite[Theorem 5.1(d) and Remark 5.2]{Auscher-Hou2024-HCL}. (Remark that $p^-_{\beta_A}(\beta)>1$ when $\beta_A<\beta<-1/2$.) \\

\paragraph{Case 2: $2<p<\infty$}
We need to prove the trace in $E^q_\delta$ and $E^{-1,q}_\delta$ for any $\delta>0$ and $q \in [p,\infty]$. The trace in $E^{-1,q}_\delta$ follows from the same arguments as in Case 1. In fact, it also works for $p=\infty$. 

Let $\beta \ge -1/2$. To prove the trace in $E^q_\delta$, by interpolation, it suffices to prove the trace in $E^p_\delta$ and $E^\infty_\delta$. To this end, we first observe that for any ball $B \subset \bR^n$, Caccioppoli's inequality (Lemma \ref{lem:Caccioppoli}) yields
\begin{equation}
    \label{e:lions-u(t)-L2(B)-Caccio}
    \begin{aligned}
        \|u(t)\|_{L^2(B)}^2
        &\lesssim t^{2\beta+1} \left( \frac{t}{r(B)^2} + 1 \right) \int_{t/2}^t \int_{2B} |s^{-(\beta+1)} u(s,y)|^2 dsdy \\
        &\quad + t^{2\beta+1} \int_{t/2}^t \int_{2B} |s^{-(\beta+1/2)} F(s,y)|^2 dsdy,
    \end{aligned}
\end{equation}
where the implicit constant does not depend on the center of $B$, nor $t$.

Let us introduce the \emph{Carleson functional} defined by
\[ \cC(u)(x) := \sup_{B:x \in B} \left( \int_0^{r(B)^2} \fint_B |u(t,y)|^2 dtdy \right)^{1/2}, \quad x \in \bR^n, \]
where $B$ describes balls in $\bR^n$. For $0<t<r(B)^2$, we infer from \eqref{e:lions-u(t)-L2(B)-Caccio} the following pointwise estimate
\begin{equation}
    \label{e:lions-u(t)-L2(B)-ptest}
    \begin{aligned}
        \|u(t)\|_{L^2(B)} 
        &\lesssim |B|^{1/2} ~ t^{\beta+1/2} \inf_{z \in B} \cC\left( \I_{(t/2,t)}(s) s^{-(\beta+1)} u(s) \right)(z) \\
        &\quad + |B|^{1/2} ~ t^{\beta+1/2} \inf_{z \in B} \cC\left( \I_{(t/2,t)}(s) s^{-(\beta+1/2)} F(s) \right)(z).
    \end{aligned}
\end{equation}

Let $\delta>0$ and $0<t<\delta$. We first prove the trace in $E^p_\delta$. For any $x \in \bR^n$, applying \eqref{e:lions-u(t)-L2(B)-ptest} to the ball $B=B(x,\delta^{1/2})$ gives
\begin{align*}
    \|u(t)\|_{E^p_\delta}^p
    &\lesssim_\delta t^{(\beta+1/2)p} \int_{\bR^n} dx \fint_{B(x,\delta^{1/2})} \left| \cC\left( \I_{(t/2,t)}(s) s^{-(\beta+1)} u(s) \right)(z) \right|^p dz \\
    &\quad + t^{(\beta+1/2)p} \int_{\bR^n} dx \fint_{B(x,\delta^{1/2})} \left| \cC\left( \I_{(t/2,t)}(s) s^{-(\beta+1/2)} F(s) \right)(z) \right|^p dz.
\end{align*}
Note that \cite[Theorem 3(b)]{Coifman-Meyer-Stein1985_TentSpaces} shows for $2<p<\infty$,
\begin{equation}
    \label{e:tent-A-C-eq}
    \|f\|_{T^p_0} \eqsim \|\cC(f)\|_{L^p(\bR^n)}. 
\end{equation}
Therefore, we obtain
\[ \|u(t)\|_{E^p_\delta} \lesssim_\delta t^{\beta+1/2} \left( \|\I_{(t/2,t)} u\|_{T^p_{\beta+1}} + \|\I_{(t/2,t)} F\|_{T^p_{\beta+1}} \right), \]
which tends to 0 as $t \to 0$ by Lebesgue's dominated convergence theorem. This proves the trace in $E^p_\delta$.

Next, consider the trace in $E^\infty_\delta$. For any ball $B \subset \bR^n$ with radius $\delta^{1/2}$, using \eqref{e:lions-u(t)-L2(B)-ptest} and \eqref{e:tent-A-C-eq} again, we obtain
\begin{align*}
    \|u(t)\|_{L^2(B)} 
    &\lesssim_\delta t^{\beta+1/2} \left( \fint_B \left| \cC\left( \I_{(t/2,t)}(s) s^{-(\beta+1)} u(s) \right)(x) \right|^p dx \right)^{1/p} \\
    &\quad + t^{\beta+1/2} \left( \fint_B \left| \cC\left( \I_{(t/2,t)}(s) s^{-(\beta+1/2)} F(s) \right)(x) \right|^p dx \right)^{1/p} \\
    &\lesssim_\delta t^{\beta+1/2} ( \|\I_{(t/2,t)} u\|_{T^p_{\beta+1}} + \|\I_{(t/2,t)} F\|_{T^p_{\beta+1}} ),
\end{align*}
which tends to 0 as $t \to 0$. This concludes the case. \\

\paragraph{Case 3: $p=\infty$}
This case is only concerned when $\beta \ge -1/2$. As aforementioned, it suffices to prove the trace in $E^{\infty}_\delta$ for $\beta>-1/2$ and in $L^2_{\loc}$ for $\beta=-1/2$. For $\beta>-1/2$, \eqref{e:lions-u(t)-L2(B)-Caccio} implies for $0<t<\delta$,
\[ \|u(t)\|_{E^\infty_\delta} \lesssim_\delta t^{\beta+1/2} \left( \|u\|_{T^\infty_{\beta+1}} + \|F\|_{T^\infty_{\beta+1/2}} \right), \]
which tends to 0 as $t \to 0$. This proves the trace in $E^\infty_\delta$.

For $\beta=-1/2$, to prove the trace in $L^2_{\loc}$, we fix a ball $B \subset \bR^n$. Let $0<t<r(B)^2$. We first claim that $u \in C([0,t];L^2(B))$. Indeed, \cite[Lemma 19]{Auscher-Hou2025_WPMRTent} implies that for $\gamma \ge 0$, $T^\infty_\gamma$ embeds into $L^2((0,T) \times B)$. As $u \in T^\infty_{1/2}$ and $\nabla u \in T^\infty_0$, we get $u \in L^2((0,t);W^{1,2}(B))$. As $\partial_t u = \Div(A\nabla u + F)$ (by \ref{item:lions_sol}), we get $\partial_t u \in L^2((0,t);W^{-1,2}(B))$. Then Lions' embedding theorem (see \textit{e.g.}, \cite[XVIII.2.1]{Dautray-JLLions1992-vol5}) yields that $u \in C([0,t];L^2(B))$ as desired. Moreover, we infer from \eqref{e:lions-u(t)-L2(B)-Caccio} that
\begin{align*}
    \sum_{j=0}^{+\infty} \|u(t/2^j)\|_{L^2(B)}^2 
    &\lesssim \int_0^t \int_{2B} |s^{-1/2} u(s,y)|^2 dsdy + \int_0^t \int_{2B} |F(s,y)|^2 dsdy \\
    &\lesssim_B \|u\|_{T^\infty_{1/2}}^2 + \|F\|_{T^\infty_0}^2.
\end{align*}
The fact that the series converges forces $\|u(t/2^j)\|_{L^2(B)}$ tends to 0 as $j \to +\infty$. By continuity, we obtain $\|u(t)\|_{L^2(B)}$ tends to 0 as $t \to 0$. 

The proof is complete.
\end{proof}

\subsection{Proof of Theorem \ref{thm:ic} \ref{item:ic-formula} and \ref{item:ic_cont}}
\label{ssec:pf-ic-cont}
Let $\beta>-1/2$ and $p_A(\beta)<p \le \infty$. Let $f \in T^p_\beta$ and $v=\cL^A_1(f)$. 

For \ref{item:ic-formula}, we first prove the identity for $f \in L^2_\ssfc(\bR^{1+n}_+)$. As $\nabla v \in T^p_{\beta+1/2}$ (by \ref{item:ic_reg}), \cite[Theorem 1.7]{Auscher-Hou2024-HCL} asserts there is a unique global weak solution to the Cauchy problem
\begin{equation}
    \label{e:ic-Delta-eq}
    \begin{cases}
        \partial_t w - \Delta w = f + \Div((A-\bI)\nabla v) \\
        w(0)=0
    \end{cases}
\end{equation}
with $\nabla w \in T^p_{\beta+1/2}$, which is given by $w=\cL^{\bI}_1(f)+\cR^{\bI}_{1/2}((A-\bI)\nabla v)$. Meanwhile, as $f \in L^2_\ssfc(\bR^{1+n}_+)$, we get from Proposition \ref{prop:formula} \ref{item:L2-ic} that $v$ itself is a weak solution to \eqref{e:ic-Delta-eq} with $\nabla v \in T^p_{\beta+1/2}$. Thus, by uniqueness, we obtain
\[ v = \cL^{\bI}_1(f)+\cR^{\bI}_{1/2}((A-\bI)\nabla v) \]
as desired in \eqref{e:LA1=LI1+RI1/2} for any $f \in L^2_\ssfc(\bR^{1+n}_+)$.

Theorem \ref{thm:ic} \ref{item:ic_reg} and Theorem \ref{thm:lions} \ref{item:lions_reg} yield all the operators involved have bounded extensions to $T^p_\beta$, so the density argument extends \eqref{e:LA1=LI1+RI1/2} to all $f \in T^p_\beta$ (or weak*-density if $p=\infty$). This proves \ref{item:ic-formula}.

To prove \ref{item:ic_cont}, we use the formula in \ref{item:ic-formula}. We infer from \cite[Theorem 1.7]{Auscher-Hou2024-HCL} and \cite[Theorem 2(e)]{Auscher-Hou2025_WPMRTent} that $\cL^\bI_1(f)$ belongs to $C([0,\infty);\scrS')$ with traces in the spaces desired in \eqref{e:ic-trace}. 

Moreover, as $\beta>-1/2$ and $\nabla v \in T^p_{\beta+1/2}$, Theorem \ref{thm:lions} \ref{item:lions_cont} yields $\cR^{\bI}_{1/2}((A-\bI)\nabla v)$ belongs to $C([0,\infty);\scrS')$ with traces also in the spaces desired in \eqref{e:ic-trace}.

Using the identity \eqref{e:LA1=LI1+RI1/2}, we obtain $v \in C([0,\infty);\scrS')$ with traces in the spaces in \eqref{e:ic-trace}. This proves \ref{item:ic_cont}. The proof is complete.
\section{Homogeneous Cauchy problem}
\label{sec:hc}

Consider the \emph{homogeneous Cauchy problem}
\begin{equation}
    \label{e:hc}
    \tag{HC}
    \begin{cases}
        \partial_t w - \Div(A\nabla w) = 0 \\
        w(0)=w_0
    \end{cases}.
\end{equation}
The main theorem of this section establishes extension of the propagator solution map $\cE_A$ defined in \eqref{e:EA} to $\DotH^{s,p}$ and hence proves the existence of weak solutions to \eqref{e:hc} asserted in Theorem \ref{thm:wp-hc}.

\begin{theorem}[Extension of $\cE_A$]
    \label{thm:hc}
    Let $\beta_A \in [-1,-1/2)$ be the constant given by Theorem \ref{thm:lions}. Let $\beta_A<\beta<0$ and
    \[ \begin{cases}
        p^-_{\beta_A}(\beta) < p \le p^+_{\beta_A}(\beta) = \infty & \text{ if } \beta \ge -1/2 \\
        p^-_{\beta_A}(\beta) < p < p^+_{\beta_A}(\beta) & \text{ if } \beta_A < \beta < -1/2 
    \end{cases}. \]
    Then $\cE_A$ extends to an operator from $\DotH^{2\beta+1,p}$ to $L^2_{\loc}((0,\infty);W^{1,2}_{\loc})$, also denoted by $\cE_A$. The following properties hold for any $w_0 \in \DotH^{2\beta+1,p}$ and $w:=\cE_A(w_0)$.
    \begin{enumerate}[label=\normalfont(\alph*)]
        \item \label{item:hc_reg} {\normalfont (Regularity)}
        $\nabla w$ belongs to $T^p_{\beta+{1}/{2}}$ with the equivalence
        \begin{equation}
            \label{e:hc_reg_nabla-w}
            \|\nabla w\|_{T^p_{\beta+{1}/{2}}} \eqsim \|w_0\|_{\DotH^{2\beta+1,p}}.
        \end{equation}
        Moreover, if $\beta_A<\beta<-1/2$, then $w$ belongs to $T^p_{\beta+1}$ with the equivalence
        \begin{equation}
            \label{e:hc_reg_w}
            \|w\|_{T^p_{\beta+1}} \eqsim \|\nabla w\|_{T^p_{\beta+{1}/{2}}} \eqsim \|w_0\|_{\DotH^{2\beta+1,p}}.
        \end{equation}

        \item \label{item:hc_formula} {\normalfont (Explicit formula)}
        It holds that
        \begin{align}
            \label{e:hc_formula_RA}
            w 
            &= \cE_{\bI}(w_0) + \cR^A_{1/2}((A-\bI) \nabla \cE_{\bI}(w_0)) \\
            \label{e:hc_formula_RDelta}
            &= \cE_{\bI}(w_0)+\cR^{\bI}_{1/2}((A-\bI)\nabla w).
        \end{align}

        \item \label{item:hc_cont} {\normalfont (Continuity)}
       	$w$ belongs to $C([0,\infty);\scrS')$ with $w(0)=w_0$.

        \item \label{item:hc_strong_cont} {\normalfont (Strong continuity)}
        For $\beta=-1/2$ and $p_-(A)<p<p_+(A)$, $w$ also belongs to $C_0([0,\infty);L^p)$ with
        \[ \sup_{t \ge 0} \|w(t)\|_{L^p} \eqsim \|w_0\|_{L^p}. \]

        \item \label{item:hc_sol}
        $w$ is a global weak solution to \eqref{e:hc} with initial data $w_0$.
    \end{enumerate}
\end{theorem}

\begin{remark}
    \label{rem:hc-w-Tp-fails}
    The equivalence \eqref{e:hc_reg_w} fails for $\beta \ge -1/2$. 
    
    For $\beta>-1/2$, as we shall see in Theorem \ref{thm:unique-u-Tpbeta+1}, any weak solution $w \in T^p_{\beta+1}$ to the equation $\partial_t w - \Div(A\nabla w) = 0$ must be null for $p_A(\beta)<p \le \infty$, even without imposing any initial condition.
    
    For $\beta=-1/2$, it also fails. Instead, the equivalence holds in a larger class, called the \emph{Kenig--Pipher space} $X^p$ introduced by \cite{Kenig-Pipher1993Xp}, which contains $T^p_{1/2}$. It has been shown in \cite[Corollaries 5.5, 5.10, and 7.5]{Auscher-Monniaux-Portal2019Lp} and \cite[Theorem 7.6]{Zaton2020wp} that the following equivalence holds
    \[ \begin{cases}
    \|w\|_{X^p} \eqsim \|\nabla w\|_{T^p_0} \eqsim \|w_0\|_{L^p} & \text{ if } p_-(A)=p^-_{\beta_A}(-1/2)<p<\infty \\ 
    \|w\|_{X^\infty} \eqsim \|w_0\|_{L^\infty},  & \text{ if } p=\infty \\
    \|\nabla w\|_{T^\infty_0} \eqsim \|w_0\|_{\bmo} & \text{ if } p=\infty 
    \end{cases}. \]
    The last equivalence corresponds to a special case $\beta=-1/2$ and $p=\infty$ in \eqref{e:hc_reg_nabla-w}.
\end{remark}

\begin{remark}
    \label{rem:hc-RA1/2}
    Let $\beta>-1$ and $\frac{n}{n+2\beta+2} < p \le \infty$. Once one can show that $\cR^A_{1/2}$ extends to a bounded operator from $T^p_{\beta+1/2}$ to $T^p_{\beta+1}$, then the proof of Theorem \ref{thm:hc} also works for initial data $w_0$ in $\DotH^{2\beta+1,p}$.
\end{remark}

To prove this theorem, we first establish the notion of distributional traces for weak solutions to the equation $\partial_t u - \Div(A\nabla u) = 0$.
\begin{prop}[Trace]
    \label{prop:trace}
    Let $\beta>-1$ and $\frac{n}{n+2\beta+2} < p \le \infty$. Let $u$ be a global weak solution to $\partial_t u - \Div(A\nabla u) = 0$ with $\nabla u \in T^p_{\beta+{1}/{2}}$. Then there exists a unique $u_{0} \in \scrS'$ so that $u(t)$ converges to $u_{0}$ in $\scrS'$ as $t \to 0$, and moreover,
    \begin{equation}
        \label{e:trace_rep}
        u = \cE_\bI(u_0)+\cR^\bI_{1/2}((A-\bI)\nabla u).
    \end{equation}
    In addition,  
    \begin{enumerate}[label=\normalfont(\roman*)]
        \item 
        \label{item:trace_beta>0}
        If $\beta \ge 0$ and $\frac{n}{n+2\beta} \le p \le \infty$, then $u_{0}$ is a constant.        
        \item 
        \label{item:trace_-1<beta<0}
        If $-1<\beta<0$, then there exist $g \in \DotH^{2\beta+1,p}$ and $c \in \bC$ so that $u_{0}=g+c$ with
        \[ \|g\|_{\DotH^{2\beta+1,p}} \lesssim \|\nabla u\|_{T^p_{\beta+{1}/{2}}}. \]
    \end{enumerate}
\end{prop}

\begin{proof}
    The proof follows from a verbatim adaptation of \cite[Proposition 6.2]{Auscher-Hou2024-HCL} for time-independent coefficient matrices. 
\end{proof}

Let us present the proof of Theorem \ref{thm:hc}.
\begin{proof}[Proof of Theorem \ref{thm:hc}]
    Let $\beta_A < \beta < 0$ and
    \[ \begin{cases}
        p^-_{\beta_A}(\beta) < p \le p^+_{\beta_A}(\beta) = \infty & \text{ if } \beta \ge -1/2 \\
        p^-_{\beta_A}(\beta) < p < p^+_{\beta_A}(\beta) & \text{ if } \beta_A < \beta < -1/2 
    \end{cases}. \]
    Define the extension of $\cE_A$ from $\DotH^{2\beta+1,p}$ to $L^2_{\loc}((0,\infty);W^{1,2}_{\loc})$ (verified below) by
    \begin{equation}
        \label{e:EA-ext}
        \cE_A(w_0) := \cE_\bI(w_0) + \cR^A_{1/2}((A-\bI) \nabla \cE_\bI(w_0)), \quad w_0 \in \DotH^{2\beta+1,p}.
    \end{equation} 
    This agrees with the formula in Proposition \ref{prop:formula} \ref{item:formula_L2_EA} when $w_0 \in L^2$. 
    
    Let us verify that for any $w_0 \in \DotH^{2\beta+1,p}$, $\cE_A(w_0)$ lies in $L^2_{\loc}((0,\infty);W^{1,2}_{\loc})$. Note that
    \begin{equation}
        \label{e:EI}
        \cE_\bI(w_0)(t,x) = (e^{t\Delta} w_0)(x), \quad (t,x) \in \bR^{1+n}_+.
    \end{equation}
    As $\DotH^{2\beta+1,p}$ embeds into $\scrS'$, we get $\cE_\bI(w_0)$ lies in $C^\infty(\bR^{1+n}_+)$, hence in $L^2_{\loc}((0,\infty);W^{1,2}_{\loc})$. Also, as $2\beta+1 < 1$, we invoke \cite[Theorem 1.1(i)]{Auscher-Hou2024-HCL} to get $\nabla \cE_\bI(w_0)$ lies in $T^p_{\beta+1/2}$ with
    \begin{equation}
        \label{e:hc-nabla-EI-w0}
        \|\nabla \cE_\bI(w_0)\|_{T^p_{\beta+1/2}} \eqsim \|w_0\|_{\DotH^{2\beta+1,p}}.
    \end{equation}
    Then applying Theorem \ref{thm:lions} \ref{item:lions_reg} gives
    \[ \|\cR^A_{1/2}((A-\bI) \nabla \cE_\bI(w_0))\|_{T^p_{\beta+1}} \lesssim \|\nabla \cE_{\bI}(w_0)\|_{T^p_{\beta+1/2}} \eqsim \|w_0\|_{\DotH^{2\beta+1,p}}, \]
    and
    \begin{equation}
        \label{e:hc-nabla-RA1/2-w0}
        \|\nabla \cR^A_{1/2}((A-\bI) \nabla \cE_\bI(w_0))\|_{T^p_{\beta+1/2}} \lesssim \|w_0\|_{\DotH^{2\beta+1,p}}.
    \end{equation}
    Thus, $\cR^A_{1/2}((A-\bI) \nabla \cE_\bI(w_0))$ lies in $L^2_{\loc}((0,\infty);W^{1,2}_{\loc})$ and so does $\cE_A(w_0)$.

    Write $w:=\cE_A(w_0)$. Let us prove the properties announced. First consider \ref{item:hc_reg}. The inequality $\|\nabla w\|_{T^p_{\beta+1/2}} \lesssim \|w_0\|_{\DotH^{2\beta+1,p}}$ directly follows from \eqref{e:hc-nabla-EI-w0} and \eqref{e:hc-nabla-RA1/2-w0}. The proof of the reverse inequality $\|w_0\|_{\DotH^{2\beta+1,p}} \lesssim \|\nabla w\|_{T^p_{\beta+1/2}}$ and \eqref{e:hc_reg_w} is deferred to the end of the proof.

    Next, consider \ref{item:hc_formula}. The first formula \eqref{e:hc_formula_RA} is exactly the definition of the extension \eqref{e:EA-ext}. For the second \eqref{e:hc_formula_RDelta}, Proposition \ref{prop:formula} \ref{item:formula_L2_EA} says it holds for $w_0 \in \DotH^{2\beta+1,p} \cap L^2$. Ensured by the \textit{a priori} estimates in \ref{item:hc_reg}, one can extend the identity (valued in $L^2_{\loc}(\bR^{1+n}_+)$) to all $w_0 \in \DotH^{2\beta+1,p}$ by density (or weak*-density if $p=\infty$). 

    For \ref{item:hc_cont}, we show the two terms in \eqref{e:EA-ext} lie in $C([0,\infty);\scrS')$. Write
    \[ v:=\cE_\bI(w_0), \quad u:=\cR^A_{1/2}((A-\bI)\nabla v). \]
    As $\nabla v \in T^p_{\beta+1/2}$, Theorem \ref{thm:lions} \ref{item:lions_cont} yields $u \in C([0,\infty);\scrS')$ with $u(0)=0$. Moreover, by \eqref{e:EI}, we get $v$ lies in $C([0,\infty);\scrS')$ with $v(0) = w_0$, and so does $w=u+v$.

    The statement \ref{item:hc_strong_cont} combines \cite[Corollary 5.10 and Proposition 5.11]{Auscher-Monniaux-Portal2019Lp} (for $p \ne 2$) and Corollary \ref{cor:Gamma} (for $p=2$).

    To prove \ref{item:hc_sol}, notice that $v$ is a global weak solution to the heat equation, and Theorem \ref{thm:lions} \ref{item:lions_sol} says $u$ is a global weak solution to $\partial_t u - \Div(A\nabla u) = \Div((A-\bI)\nabla v)$. Thus, $w=u+v$ is a global weak solution to $\partial_t w - \Div(A\nabla w) = 0$. Moreover, \ref{item:hc_cont} yields $w(t)$ converges to $w_0$ as $t \to 0$ in $\scrS'$, hence in $\scrD'$. This proves \ref{item:hc_sol}.

    Let us finish by proving the rest of \ref{item:hc_reg}. First, we show the reverse inequality $\|w_0\|_{\DotH^{2\beta+1,p}} \lesssim \|\nabla w\|_{T^p_{\beta+1/2}}$ in \eqref{e:hc_reg_nabla-w}. Note that $w$ is a global weak solution to $\partial_t w - \Div(A\nabla w) = 0$ with $\nabla w \in T^p_{\beta+1/2}$, and $w(t)$ converges to $w_0$ in $\scrS'$ as $t \to 0$. Hence, $w_0$ agrees with the trace $g$ of $w$ constructed in Proposition \ref{prop:trace}, which gives $\|w_0\|_{\DotH^{2\beta+1,p}} \lesssim \|\nabla w\|_{T^p_{\beta+1/2}}$. 

    It only remains to prove \eqref{e:hc_reg_w}. Using again the fact that $w$ is a global weak solution to $\partial_t w - \Div(A\nabla w) = 0$, we deduce from Corollary \ref{cor:Caccioppoli-tent} that $\|\nabla w\|_{T^p_{\beta+1/2}} \lesssim \|w\|_{T^p_{\beta+1}}$. 
    
    On the other hand, \cite[Theorem 3.2]{Auscher-Hou2024-HCL} shows that when $\beta<-1/2$ (\textit{i.e.}, $2\beta+1<0$), $\cE_\bI(w_0)$ belongs to $T^p_{\beta+1}$ with the equivalence
    \[ \|\cE_\bI(w_0)\|_{T^p_{\beta+1}} \eqsim \|w_0\|_{\DotH^{2\beta+1,p}}. \]
    Then we infer from \eqref{e:EA-ext} that
    \begin{align*}
        \|w\|_{T^p_{\beta+1}} 
        &\le \|\cE_\bI(w_0)\|_{T^p_{\beta+1}} + \|\cR^A_{1/2}((A-\bI) \nabla \cE_\bI(w_0))\|_{T^p_{\beta+1}} \\
        &\lesssim \|w_0\|_{\DotH^{2\beta+1,p}} + \|\nabla \cE_\bI(w_0)\|_{T^p_{\beta+1/2}} \lesssim \|w_0\|_{\DotH^{2\beta+1,p}}.
    \end{align*}
    This proves \ref{item:hc_reg}. The proof is complete.
\end{proof}

\begin{cor}[Continuity of propagators]
    \label{cor:Gamma-cont}
    Let $p_-(A)<p<p_+(A)$ and $p_-(A)'<q<p_+(A)'$. Let $g \in L^p$ and $h \in L^q$.
    \begin{enumerate}[label=\normalfont(\roman*)]
        \item \label{item:t-Gamma-cont}
        The function $t \mapsto \Gamma_A(t,s)g$ lies in $C_0([s,\infty);L^p)$.
        
        \item \label{item:s-Gamma-w-cont}
        The function $s \mapsto \Gamma_A(t,s)g$ lies in $C_\ssfw([0,t];L^p)$. \footnote{Here, $C_\ssfw([0,t];E)$ is the space of continuous functions valued in $E$ equipped with its weak topology.}

        \item \label{item:t-Gamma*-w-cont}
        The function $t \mapsto \Gamma_A(t,s)^\ast h$ lies in $C_\ssfw([s,\infty);L^q)$.

        \item \label{item:s-Gamma*-cont}
        The function $s \mapsto \Gamma_A(t,s)^\ast h$ lies in $C([0,t];L^q)$.
    \end{enumerate}
\end{cor}

\begin{proof}
    By duality, it suffices to prove \ref{item:t-Gamma-cont} and \ref{item:s-Gamma*-cont}. The statement \ref{item:t-Gamma-cont} follows from Theorem \ref{thm:hc} \ref{item:hc_strong_cont} by shifting the time. We present two approaches to prove \ref{item:s-Gamma*-cont}. Both of them use the duality relation \eqref{e:Gamma-dual} in Corollary \ref{cor:Gamma-}, which says
    \[ \Gamma_A(t,s)^\ast = \Gamma^-_{A^\ast}(s,t) = \Gamma_{\tilA_t}(t-s,0), \quad 0 \le s \le t, \]
    where $\tilA_t$ is defined in \eqref{e:tilAt} given by
    \[ \tilA_t(s,x) := 
        \begin{cases}
            A^\ast(t-s,x), & \text{ if } 0 \le s \le t \\
            \Lambda_0 \bI & \text{ if } s>t
        \end{cases}. \]
    
    The first way is to adapt the proof of \ref{item:t-Gamma-cont} to the backward equation. Detailed computation is left to the reader.

    The second proof is based on the observation that
    \begin{equation}
        \label{e:p-+(tilAt)-p-+(A)'}
        p_-(\tilA_t) \le p_+(A)'<p_-(A)' \le p_+(\tilA_t).
    \end{equation}
    Indeed, since $(\Gamma_A(t,s))$ is uniformly bounded on $L^p$ for $p_-(A)<p<p_+(A)$, by duality, $(\Gamma_A(t,s)^\ast)$ is uniformly bounded on $L^q$ for $p_+(A)'<q<p_-(A)'$, and so is $(\Gamma_{\tilA_t}(\tau,0))_{0 \le \tau \le t}$. Moreover, as $p_-(\Lambda_0 \bI)=1$ and $p_+(\Lambda_0 \bI)=\infty$, we get $(\Gamma_{\tilA_t}(\tau,0))_{\tau \ge 0}$ is uniformly bounded on $L^q$ for $p_+(A)'<q<p_-(A)'$. Hence, \eqref{e:p-+(tilAt)-p-+(A)'} follows by shifting the time.

    Theorem \ref{thm:hc} \ref{item:hc_strong_cont} yields for $p_-(\tilA_t)<q<p_+(\tilA_t)$ (in particular for $p_+(A)'<q<p_-(A)'$) and $\psi \in L^q$, the function $\tau \mapsto \Gamma_{\tilA_t}(\tau,0)\psi$ belongs to $C_0([0,\infty);L^q)$. Applying it to the function $s \mapsto \Gamma_{\tilA_t}(t-s,0) h = \Gamma_A(t,s)^\ast h$ for $0 \le s \le t$ hence gives \ref{item:s-Gamma*-cont}. This completes the proof.
\end{proof}
\section{Uniqueness and representation}
\label{sec:unique}

Let us first state two theorems on uniqueness of weak solutions. The proof of the representation theorem (cf. Theorem \ref{thm:rep}) is presented in Section \ref{ssec:pf-rep}.

\begin{theorem}[Uniqueness]
    \label{thm:unique}
    Let $\beta>\beta_A$ and
    \[ \begin{cases}
        p^-_{\beta_A}(\beta) < p \le p^+_{\beta_A}(\beta) = \infty & \text{ if } \beta \ge -1/2 \\
        p^-_{\beta_A}(\beta) < p < p^+_{\beta_A}(\beta) & \text{ if } \beta_A < \beta < -1/2 
    \end{cases}. \]
    Let $u$ be a global weak solution to the Cauchy problem
    \[ \begin{cases}
        \partial_t u-\Div(A\nabla u) = 0, \\ 
        u(0)=0
    \end{cases}, \]
    with $\nabla u \in T^p_{\beta+{1}/{2}}$. Then $u=0$.
\end{theorem}

We also present another class for which the uniqueness holds without imposing any initial condition.
\begin{theorem}
    \label{thm:unique-u-Tpbeta+1}
    Let $\beta>-1/2$ and $p_A(\beta) < p \le \infty$. Let $u \in T^p_{\beta+1}$ be a global weak solution to the equation $\partial_t u - \Div(A\nabla u) = 0$. Then $u=0$.
\end{theorem}
The proofs are presented right below in Sections \ref{ssec:pf-unique} and \ref{ssec:pf-unique-Tpbeta+1}.

\subsection{Proof of Theorem \ref{thm:unique}}
\label{ssec:pf-unique}
The following lemma reduces the proof to the case $p>p_-(A)$. Note that $p \le p_-(A)$ occurs only if $\beta>-1/2$.

\begin{lemma}
    \label{lem:p<p_(A)-embed}
    Let $\beta>-1/2$ and $p_A(\beta)<p \le p_-(A)$. Then there exist $\beta_0>-1/2$ and $p_0 \in (p_-(A),2)$ so that $T^p_{\beta+\gamma}$ embeds into $T^{p_0}_{\beta_0+\gamma}$ for any $\gamma \in \bR$.
\end{lemma}

\begin{proof}
    By embedding of tent spaces, it suffices to find $\beta_0 \in (-1/2,\beta)$ and $p_0 \in (p_-(A),2)$ so that $2\beta_0-\frac{n}{p_0} = 2\beta-\frac{n}{p}$. To do so, first pick $\beta_1$ with $2\beta_1-\frac{n}{p_-(A)} = 2\beta-\frac{n}{p}$. As $p \le p_-(A)$, we have $\beta_1 \le \beta$. We claim $\beta_1>-1/2$. If it holds, then perturbation gives $\beta_0 \in (-1/2,\beta_1) \subset (-1/2,\beta)$ and $p_0 \in (p_-(A),2)$ with $2\beta_0-\frac{n}{p_0} = 2\beta-\frac{n}{p}$.
    
    It only remains to verify the claim. Note that 
    \[ r>p_A(\gamma)=\frac{np_-(A)}{n+(2\gamma+1)p_-(A)} \iff 2\gamma-\frac{n}{r}>-1-\frac{n}{p_-(A)}. \] 
    As $p>p_A(\beta)$, we get $2\beta_1 = 2\beta-\frac{n}{p} + \frac{n}{p_-(A)} > -1$, so $\beta_1>-1/2$.
\end{proof}

Let us present the proof of Theorem \ref{thm:unique}. Let $u$ be a global weak solution to $\partial_t u - \Div(A\nabla u) = 0$ with $\nabla u \in T^p_{\beta+1/2}$ and null initial data. The argument is split into 2 cases: $\beta \ge -1/2$ and $\beta<-1/2$.

\subsubsection{\texorpdfstring{$\beta \ge -1/2$}{beta>=-1/2}}
\label{sssec:unique-beta>=-1/2}
Thanks to Lemma \ref{lem:p<p_(A)-embed}, it suffices to consider the case $p_-(A)<p \le \infty$. We infer from \cite[Lemma 3.7]{Auscher-Hou2024-HCL} that $u$ satisfies the local $L^2$-norm estimate that for $0<a<b<\infty$ and $R>1$,
\[ \int_a^b \int_{B(0,R)} |u|^2 \lesssim_{a,b,p,\beta} R^{3n+2} \left( \|\nabla u\|_{T^p_{\beta+1/2}}^2 + \|u\|_{L^2((a,b) \times B(0,1))}^2 \right). \]
In fact, it holds for any $u \in L^2_{\loc}(\bR^{1+n}_+)$ with $\nabla u \in T^p_{\beta+1/2}$. This allows us to invoke \cite[Theorem 5.1]{Auscher-Monniaux-Portal2019Lp} to get an interior representation of $u$, also called ``\emph{homotopy identity}", that for $0<s<t<\infty$ and $h \in \Cc(\bR^n)$,
\begin{equation}
    \label{e:homo-id}
    \int_{\bR^n} u(t,x) \ovh(x) dx = \int_{\bR^n} u(s,x) \overline{\Gamma_A(t,s)^\ast h}(x) dx.
\end{equation}

Meanwhile, as $p>p_A(\beta)>\frac{n}{n+2\beta+2}$ and $u(0)=0$, Proposition \ref{prop:trace} yields
\begin{equation}
    \label{e:unique-u=RI1/2}
    u = \cR^\bI_{1/2}((A-\bI)\nabla u).
\end{equation}
Hence, Theorem \ref{thm:lions} \ref{item:lions_cont} implies $u(s)$ tends to 0 as $s \to 0$ in $\scrS'$ and in some finer trace spaces depending on $\beta$ and $p$. 

Therefore, our strategy is to analyze the limit $s \to 0$ for \eqref{e:homo-id}, by employing these trace spaces and the continuity of the propagators, to obtain $u(t)=0$ in $\scrD'$ for any $t>0$, and hence $u=0$. The argument is split into 3 sub-cases. \\

\paragraph{Case 1: $p_-(A)<p \le 2$}
Using \eqref{e:unique-u=RI1/2} and Theorem \ref{thm:lions} \ref{item:lions_cont}, we get
\[ \lim_{s \to 0} u(s) = 0 \quad \text{ in } L^p. \]
Moreover, Corollary \ref{cor:Gamma-cont} \ref{item:s-Gamma*-cont} says
\[ \lim_{s \to 0} \Gamma_A(t,s)^\ast h = \Gamma_A(t,0)^\ast h \quad \text{ in } L^{p'}. \]
Then taking the limit $s \to 0$ in \eqref{e:homo-id} gives $u(t)=0$ as desired. \\

\paragraph{Case 2: $2<p<\infty$}
Let $\delta>0$. Using \eqref{e:unique-u=RI1/2} and Theorem \ref{thm:lions} \ref{item:lions_cont} again, we get
\[ \lim_{s \to 0} u(s) = 0 \quad \text{ in } E^\infty_\delta. \]
For the other term, we invoke \cite[Lemma 4.7(2)]{Auscher-Monniaux-Portal2019Lp} to get
\[ \lim_{s \to 0} \Gamma_A(t,s)^\ast h = \Gamma_A(t,0)^\ast h \quad \text{ in } E^1_\delta. \]  
By duality, taking the limit $s \to 0$ in \eqref{e:homo-id} gives $u(t)=0$ as desired. \\

\paragraph{Case 3: $p=\infty$}
For $\beta>-1/2$, it follows exactly the same arguments as in Case 2. To prove the case $\beta=-1/2$, we also use \eqref{e:homo-id}. Let $h \in \Cc(\bR^n)$. Pick $M>1$ so that $\supp(h) \subset B(0,M) =: B$. We claim that there exists a constant $c>0$ so that for any $x \in \bR^n$,
\begin{equation}
    \label{e:unique-Gamma*h-pt-bd}
    \sup_{0 \le s \le t/2} \|\I_{B(x,1)} \Gamma_A(t,s)^\ast h\|_{L^2} \lesssim \|h\|_{L^2} \left( \I_{4B}(x) + e^{-c\frac{|x|^2}{t}} \I_{(4B)^c}(x) \right).
\end{equation}
Indeed, when $x \in 4B$, using $L^2$-boundedness of $(\Gamma_A(t,s)^\ast)$, we have
\[ \sup_{0 \le s \le t/2} \|\I_{B(x,1)} \Gamma_A(t,s)^\ast h\|_{L^2} \le \sup_{0 \le s \le t/2} \|\Gamma_A(t,s)^\ast h\|_{L^2} \le \|h\|_{L^2}. \]
When $x \notin 4B$, we get $\dist(B(x,1),B) \ge \frac{1}{2}|x|$. Using the $L^2-L^2$ off-diagonal decay of $(\Gamma_A(t,s)^\ast)$, we obtain a constant $c>0$ so that
\[ \|\I_{B(x,1)} \Gamma_A(t,s)^\ast h\|_{L^2} \lesssim e^{-4c\frac{\dist(B(x,1),B)^2}{t-s}} \|h\|_{L^2} \le e^{-c\frac{|x|^2}{t}} \|h\|_{L^2} \]
as desired. This proves \eqref{e:unique-Gamma*h-pt-bd}.

Then taking averages over $B(x,1)$ for the integral on the right-hand side of \eqref{e:homo-id}, we obtain
\begin{equation}
    \label{e:unique-endpt-av}
    \int_{\bR^n} u(t,x) \ovh(x) dx = \int_{\bR^n} dx \fint_{B(x,1)} u(s,y) \overline{\Gamma_A(t,s)^\ast h} (y) dy.
\end{equation}
Remark \ref{rem:lions-endpt-Einfty-bd} says $u(s)$ is uniformly bounded in $E^\infty_1$ for $0<s<1$. Combining it with \eqref{e:unique-Gamma*h-pt-bd} gives that for $0<s<\min\{t/2,1\}$ and $x \in \bR^n$,
\[ \fint_{B(x,1)} |u(s)| |\Gamma_A(t,s)^\ast h| \lesssim \|u\|_{E^\infty_1} \|h\|_{L^2} \left( \I_{4B}(x) + e^{-c\frac{|x|^2}{t}} \I_{(4B)^c}(x) \right), \]
which is integrable on $\bR^n$. Moreover, Theorem \ref{thm:lions} \ref{item:lions_cont} says $u(s)$ tends to 0 as $s \to 0$ in $L^2_{\loc}$, and Corollary \ref{cor:Gamma-cont} \ref{item:s-Gamma*-cont} says $\Gamma_A(t,s)^\ast h$ converges to $\Gamma_A(t,0)^\ast h$ as $s \to 0$ in $L^2$. Therefore, we obtain that for any $x \in \bR^n$,
\[ \lim_{s \to 0} \fint_{B(x,1)} u(s,y) \overline{\Gamma_A(t,s)^\ast h}(y) dy = 0. \]
Applying Lebesgue's dominated convergence theorem to the limit $s \to 0$ of the integral on the right-hand side of \eqref{e:unique-endpt-av} implies $u(t)=0$. 
    
This concludes the case $\beta \ge -1/2$.

\subsubsection{\texorpdfstring{$\beta<-1/2$}{beta<-1/2}}
\label{sssec:unique-beta<-1/2}
In this case, we use the interpolation argument as in the proof of Lemma \ref{lem:lions-beta<-1/2-p=2}. Recall the map $\Phi_A$ defined in \eqref{e:PhiA} given by
\begin{align*}
    \Phi_A : S^p_{\beta+1/2} &\to \Div T^p_{\beta+1/2} \times (\DotH^{2\beta+1,p}+\bC) \\
    \Phi_A(u) &:= (\partial_t u -\Div(A\nabla u), \quad \tr(u)).
\end{align*}
Lemma \ref{lem:tr-Sp} implies $\Phi_A$ is bounded for $-1<\beta<0$ and $1<p<\infty$. Moreover, we have also shown that $\Phi_A$ is an isomorphism for $\beta_A<\beta \le -1/2$ and $p=2$, see Section \ref{ssec:pf-lions-beta<-1/2-p=2}. Let us recall the construction of the inverse $\Psi_A$ defined in \eqref{e:PsiA(g,v0)}. For any $g \in \Div T^p_{\beta+1/2}$ and $v_0 \in \DotH^{2\beta+1,p}+\bC$, $v=:\Psi_A(g,v_0)$ is the unique weak solution in $S^p_{\beta+1/2}$ to the Cauchy problem
\begin{equation}
    \label{e:Psi(g,v0)-eq-unique}
    \begin{cases}
        \partial_t v - \Div(A\nabla v) = g, \\
        v(0)=v_0
    \end{cases}.
\end{equation}

Let $\beta=-1/2$ and $p_-(A) = p^-_{\beta_A}(-1/2)<p<p^+_{\beta_A}(-1/2) = \infty $. For any $g \in \Div T^p_{0}$ and $v_0 \in \DotH^{0,p}+\bC$, Theorems \ref{thm:lions} and \ref{thm:hc} yield there is a weak solution $v \in S^p_{0}$ to \eqref{e:Psi(g,v0)-eq-unique} with
\[ \|v\|_{S^p_{0}} \lesssim \|g\|_{\Div T^p_{0}} + \|v_0\|_{\DotH^{0,p}}. \]
The proof in Section \ref{sssec:unique-beta>=-1/2} ensures the weak solution is unique in $S^p_{0}$. Therefore, $\Psi_A$ is bounded for $\beta=-1/2$ and $p_-(A)<p<\infty$. 

Then by interpolation, $\Psi_A$ extends to a bounded map from $\Div T^p_{\beta+1/2} \times (\DotH^{2\beta+1,p}+\bC)$ to $S^p_{\beta+1/2}$ for $\beta_A<\beta \le -1/2$ and $p^-_{\beta_A}(\beta) < p < p^+_{\beta_A}(\beta)$. 
\footnote{In fact, since both spaces are semi-normed, one needs to pass to the quotient map (from $\Div T^p_{\beta+1/2} \times \DotH^{2\beta+1,p}/\bC$ to $S^p_{\beta+1/2}/\bC$) to use interpolation, and then lifts up, as shown in Figure \ref{fig:Phi_A-diagram}. We leave the detailed verification to the reader.}
By density, one also gets $\Psi_A$ is the inverse of $\Phi_A$. In particular, for $g=0$ and $v_0=0$, it implies that the unique weak solution $v \in S^p_{\beta+1/2}$ to the Cauchy problem
\begin{equation}
    \label{e:unique-eq}
    \begin{cases}
        \partial_t v - \Div(A\nabla v) = 0, \\
        v(0)=0
    \end{cases}
\end{equation}
must be null. Observe that any global weak solution $v$ to \eqref{e:unique-eq} with $\nabla v \in T^p_{\beta+1/2}$ belongs to $S^p_{\beta+1/2}$, using the equation. We hence obtain the uniqueness in the class $\nabla v \in T^p_{\beta+1/2}$. This completes the proof.

\subsection{Proof of Theorem \ref{thm:unique-u-Tpbeta+1}}
\label{ssec:pf-unique-Tpbeta+1}
Let $\beta>-1/2$ and $p_A(\beta)<p \le \infty$. Let $u \in T^p_{\beta+1}$ be a global weak solution to the equation $\partial_t u-\Div(A\nabla u)=0$. Using Lemma \ref{lem:p<p_(A)-embed} again, we may assume $p_-(A)<p \le \infty$. Still, our strategy is to use the homotopy identity \eqref{e:homo-id}. Note that as $u \in T^p_{\beta+1}$, \cite[Lemmas 13 and 19]{Auscher-Hou2025_WPMRTent} yield that for $0<a<b<\infty$ and $R \ge 1$,
\[ \int_a^b \int_{B(0,R)} |u|^2 \lesssim_{a,b,p,\beta} R^n \left( \|u\|_{T^p_{\beta+1}}^2 + \|u\|_{L^2((a,b) \times B(0,b^{1/2})}^2 \right), \]
which verifies the conditions in \cite[Theorem 5.1]{Auscher-Monniaux-Portal2019Lp} to obtain the homotopy identity \eqref{e:homo-id}. 

It remains to take the limits $s \to 0$ in \eqref{e:homo-id} to get $u(t)=0$ in $\scrD'$ for any $t>0$ and hence $u=0$. To show this, we first observe that for any $s>0$, $u(s)$ lies in $E^p_{s/16}$ with
\begin{equation}
    \label{e:unique-Tpbeta+1-u(s)Eps/16}
    \|u(s)\|_{E^p_{s/16}} \lesssim s^{\beta+1/2} \|u\|_{T^p_{\beta+1}}.
\end{equation}
Indeed, for $p<\infty$, it follows by applying Caccioppoli's inequality (cf. Lemma \ref{lem:Caccioppoli}) to the average of $u(s)$ on $B(x,\frac{\sqrt{s}}{4})$ as 
\[ \|u(s)\|_{E^p_{s/16}} \lesssim \left( \int_{\bR^n} \left( \frac{1}{s} \int_{s/2}^s \fint_{B(x,\frac{\sqrt{s}}{2})} |u|^2 \right)^{p/2} dx \right)^{1/p} \lesssim s^{\beta+1/2} \|u\|_{T^p_{\beta+1}}. \]
The same argument also works for $p=\infty$. Then the discussion is split into two cases. \\

\paragraph{Case 1: $p_-(A)<p \le 2$}
H\"older's inequality yields that $E^p_{s/16}$ embeds into $L^p$, so we get
\[ \|u(s)\|_{L^p} \lesssim \|u(s)\|_{E^p_{s/16}} \lesssim s^{\beta+1/2} \|u\|_{T^p_{\beta+1}}, \]
which tends to 0 as $s \to 0$. Meanwhile, Corollary \ref{cor:Gamma-cont} \ref{item:s-Gamma*-cont} says
\[ \lim_{s \to 0} \Gamma_A(t,s)^\ast h = \Gamma(t,0)^\ast h \quad \text{ in } L^{p'}. \]
Then taking the limit $s \to 0$ in \eqref{e:homo-id} implies $u(t)=0$ as desired. \\

\paragraph{Case 2: $2<p \le \infty$}
Let $\delta>0$ and $0<s<t<\delta$. Using the change of aperture, cf. \eqref{e:Ep-aperture}, we get $E^p_{s/16}$ embeds into $E^p_\delta$ with
\[ \|u(s)\|_{E^p_\delta} \lesssim \|u(s)\|_{E^p_{s/16}} \lesssim s^{\beta+1/2} \|u\|_{T^p_{\beta+1}}, \]
which tends to 0 as $s \to 0$. On the other hand, \cite[Lemma 4.7(2)]{Auscher-Monniaux-Portal2019Lp} says
\[ \lim_{s \to 0} \Gamma_A(t,s)^\ast h = \Gamma_A(t,0)^\ast h \quad \text{ in } E^{p'}_\delta. \]
Taking the limit $s \to 0$ in \eqref{e:homo-id} implies $u(t)=0$. 

This completes the proof.

\subsection{Proof of Theorem \ref{thm:rep}}
\label{ssec:pf-rep}
Let $\beta>\beta_A$ and 
\[ \begin{cases}
    p^-_{\beta_A}(\beta) < p \le p^+_{\beta_A}(\beta) = \infty & \text{ if } \beta \ge -1/2 \\
    p^-_{\beta_A}(\beta) < p < p^+_{\beta_A}(\beta) & \text{ if } \beta_A < \beta < -1/2 
\end{cases}. \]
Let $u$ be a global weak solution to $\partial_t u - \Div(A\nabla u) = 0$ with $\nabla u \in T^p_{\beta+1/2}$. As $p>p^-_{\beta_A}(\beta)>\frac{n}{n+2\beta+2}$, Proposition \ref{prop:trace} yields there exists a unique $u_0 \in \scrS'$ so that $u(t)$ converges to $u_0$ as $t \to 0$, and
\[ u = \cE_\bI(u_0) + \cR^\bI_{1/2}((A-\bI)\nabla u). \]

Then we prove the properties of the trace. First consider \ref{item:rep_beta>0}. For $\beta \ge 0$ and $\frac{n}{n+2\beta} \le p \le \infty$, Proposition \ref{prop:trace} \ref{item:trace_beta>0} says $u_0$ is a constant, and so is $\cE_\bI(u_0)$. Thus, $v:=u-\cE_\bI(u_0)$ is a global weak solution to the Cauchy problem
\[ \begin{cases}
    \partial_t v - \Div(A\nabla v) = 0 \\
    v(0) = 0
\end{cases} \]
with $\nabla v \in T^p_{\beta+1/2}$. As $p^-_{\beta_A}(\beta)<p \le p^+_{\beta_A}(\beta)=\infty$, Theorem \ref{thm:unique} yields $v=0$, so $u=\cE_\bI(u_0)$ is a constant. This proves \ref{item:rep_beta>0}.

Next, to prove \ref{item:rep_-1<beta<0}, note that Proposition \ref{prop:trace} \ref{item:trace_-1<beta<0} says there exist $g \in \DotH^{2\beta+1,p}$ and $c \in \bC$ so that $u_0=g+c$. Then we use the formula \eqref{e:hc_formula_RDelta} proved in Theorem \ref{thm:hc} \ref{item:hc_formula} to get
\[ u = \cE_\bI(g) + \cR^\bI_{1/2}((A-\bI)\nabla u) + c = \cE_A(g)+c. \]
This proves \ref{item:rep_-1<beta<0}. The proof is complete.
\begin{remark}
    \label{rem:rep-system-ellp}
    The same proof applies to parabolic systems under the Gårding ellipticity condition \eqref{e:Garding-elliptic}. In particular, for $\beta=-1/2$ ($s=0$) and $p_-(A)<p<2$, it reduces the pointwise (strong) ellipticity condition used in \cite[Theorem 1.1]{Zaton2020wp}.
\end{remark}
\section{Results for homogeneous Besov spaces}
\label{sec:Besov}

This section discusses extension of our results to initial data taken in homogeneous Besov spaces $\DotB^s_{p,p}$. The definition of $\DotB^s_{p,p}$ can be analogously adapted from Definition \ref{def:Hsp}, as ``realizations'' of homogeneous Besov spaces defined on $\scrS'/\scrP$. The counterparts of tent spaces are $Z$-spaces (introduced by \cite{Barton-Mayboroda2016_Zspaces} with a different notion). For any $\beta \in \bR$ and $p \in (0,\infty)$, the \emph{(parabolic) $Z$-space} $Z^p_\beta$ consists of measurable functions $f$ on $\bR^{1+n}_+$ for which
\[ \|f\|_{Z^p_\beta} := \left( \int_0^\infty \int_{\bR^n} \left( \int_{t/2}^t \fint_{B(x,t^{1/2})} |s^{-\beta} f(s,y)|^2 dsdy \right)^{p/2} \frac{dt}{t} dx \right)^{1/p} < \infty. \]
For $p=\infty$, the norm is given by
\[ \|f\|_{Z^{\infty}_\beta} := \sup_{t>0, ~ x \in \bR^n} \left( \int_{t/2}^t \fint_{B(x,t^{1/2})} |s^{-\beta} f(s,y)|^2 dsdy \right)^{1/2}. \]

The relation between homogeneous Besov spaces (resp. $Z$-spaces) and homogeneous Hardy--Sobolev spaces (resp. tent spaces) are given by real interpolation. Namely, let $0<p_0,p_1 \le \infty$ and $s_0,s_1,\beta_0,\beta_1 \in \bR$ with $s_0 \ne s_1$ and $\beta_0 \ne \beta_1$. Suppose $s=(1-\theta) s_0 + \theta s_1$, $\beta = (1-\theta) \beta_0 + \theta \beta_1$, and $\frac{1}{p} = \frac{1-\theta}{p_0} + \frac{\theta}{p_1}$ for some $\theta \in (0,1)$. Then one gets
\begin{equation}
    \label{e:real-inter-B-Z}
    (\DotH^{s_0,p_0}, \DotH^{s_1,p_1})_{\theta,p} = \DotB^s_{p,p}, \quad (T^{p_0}_{\beta_0}, T^{p_1}_{\beta_1})_{\theta,p} = Z^p_\beta.
\end{equation}
See \cite[Theorem 2.30]{Amenta-Auscher2018-book}. 

We also have the Sobolev embedding, see \textit{e.g.}, \cite[Theorem 2.34]{Amenta-Auscher2018-book}. Let $0<p_0<p_1<p_2 \le \infty$ and $\beta_0>\beta_1>\beta_2$ with $2\beta_0-\frac{n}{p_0} = 2\beta_1-\frac{n}{p_1} = 2\beta_2-\frac{n}{p_2}$. Then
\begin{equation}
    \label{e:sobolev-embed-T-Z}
    T^{p_0}_{\beta_0} \hookrightarrow Z^{p_1}_{\beta_1} \hookrightarrow T^{p_2}_{\beta_2}.
\end{equation}

Using these tools, \cite[Theorem 9.2]{Auscher-Hou2024-HCL} shows that for $-1<s<1$ and $\frac{n}{n+s+1} \le p \le \infty$, the solution map $u_0 \mapsto u(t,x) = (e^{t\Delta} u_0)(x)$ is an isomorphism from $\DotB^s_{p,p} + \bC$ to the space of distributional solutions $u$ to the heat equation (with null source terms) with $\nabla u \in Z^p_{s/2}$. Particularly, the equivalence holds that
\[ \|\nabla e^{t\Delta} u_0 \|_{Z^p_{s/2}} \eqsim \|u_0\|_{\DotB^s_{p,p}}. \]

Moreover, at the extreme point $s=1$ and $p=\infty$, \cite[Proposition 9.1]{Auscher-Hou2024-HCL} shows that the solution map is an isomorphism from the \emph{homogeneous Sobolev space} $\DotW^{1,\infty}$ (rather than $\DotB^1_{\infty,\infty}$) to the space of distributional solutions $u$ to the heat equation with $\nabla u \in Z^\infty_{1/2}$. The equivalence also holds that
\begin{equation}
    \label{e:heat-besov-eq-s=1-p=infty}
    \|\nabla e^{t\Delta} u_0\|_{Z^\infty_{1/2}} \eqsim \|u_0\|_{\DotW^{1,\infty}}.
\end{equation}
We call this change as the \emph{canonical modification} for $s=1$ ($\beta=0$) and $p=\infty$. Recall that $\DotW^{1,\infty}$ consists of distributions $\varphi$ for which $\nabla \varphi$ is bounded, endowed with the semi-norm $\|\varphi\|_{\DotW^{1,\infty}} := \|\nabla \varphi\|_{L^\infty}$. Thanks to Rademacher's theorem, we know that $\DotW^{1,\infty}$ agrees with the set of Lipschitz continuous functions, up to almost everywhere equality.

Analogous results are also established there for parabolic Cauchy problems of type \eqref{e:NaPC} with \emph{time-independent} coefficients.

The following theorem extends the above results to the non-autonomous case.
\begin{theorem}[Well-posedness of Cauchy problems of type \eqref{e:NaPC} for homogeneous Besov spaces]
    \label{thm:besov}
    Let $\beta_A \in [-1,-1/2)$ be the constant given by Theorem \ref{thm:lions}. Let $\beta>\beta_A$ and
    \[ \begin{cases}
        p^-_{\beta_A}(\beta) < p \le p^+_{\beta_A}(\beta) = \infty & \text{ if } \beta > -1/2 \\
        p^-_{\beta_A}(\beta) < p < p^+_{\beta_A}(\beta) & \text{ if } \beta_A < \beta \le -1/2 
    \end{cases}. \]
    With the canonical modification for $\beta=0$ and $p=\infty$ (\textit{i.e.}, the initial data $u_0$ is taken in $\DotW^{1,\infty}$), Theorems \ref{thm:ic}, \ref{thm:lions}, \ref{thm:hc} (except \ref{item:hc_strong_cont}), \ref{thm:unique}, \ref{thm:unique-u-Tpbeta+1}, and \ref{thm:rep} are all valid when replacing homogeneous Hardy--Sobolev spaces (resp. tent spaces) by homogeneous Besov spaces (resp. $Z$-spaces) with the same exponents.
\end{theorem}

We exclude here the extreme point $\beta=-1/2$ and $p=\infty$, since we do not know whether the Lions operator $\cR^A_{1/2}$ is bounded from $Z^\infty_0$ to $Z^\infty_{1/2}$. Moreover, it is not clear whether weak solutions are continuous in $\DotB^0_{p,p}$ (except in the case $p=2$ as $\DotB^0_{2,2}$ identifies with $L^2$), which would be the analog of Theorem \ref{thm:hc} \ref{item:hc_strong_cont}.

\begin{proof} 
	We just show some main ingredients of the proof. For convenience, we still use the same label of the theorems for their variants in Besov spaces. 

    First consider Theorem \ref{thm:lions}. It follows from real interpolation \eqref{e:real-inter-B-Z} and properties of the extension of $\cR^{\bI}_{1/2}$ on $Z$-spaces, see \cite[Theorem 9.2]{Auscher-Hou2024-HCL}. We only need to prove the trace space $E^q_\delta$ in \ref{item:lions_cont} for $\beta=-1/2$ and $2<p<\infty$. Let $F \in Z^p_0$ and $u:=\cR^A_{1/2}(F)$. Using interpolation of slice spaces, it suffices to show $u(t)$ tends to 0 as $t \to 0$ in $E^p_\delta$ and $E^\infty_\delta$ for any $\delta>0$. Let $0<t<\delta$. To prove the trace in $E^p_\delta$, we use change of aperture of slice spaces (cf. \eqref{e:Ep-aperture}) and Caccioppoli's inequality (cf. Lemma \ref{lem:Caccioppoli}) to get
    \begin{align*}
        \|u(t)\|_{E^p_\delta}^p 
        &\lesssim \|u(t)\|_{E^p_{t/16}}^p \\
        &\lesssim \int_{\bR^n} \left( \int_{t/2}^t \fint_{B(x,\sqrt{t}/2)} |s^{-1/2} u(s,y)|^2 dsdy \right)^{p/2} dx =: \Phi(t)
    \end{align*}
    Notice that
    \[ \int_0^\infty \Phi(t) \frac{dt}{t} \lesssim \|u\|_{Z^p_{1/2}}^p < \infty. \]
    The fact that the integral converges forces that $\Phi(t)$ tends to 0 as $t \to 0$. This proves the trace in $E^p_\delta$.

    To prove the trace in $E^\infty_\delta$, we use Lions' embedding theorem. Let $B \subset \bR^n$ be a ball with radius $\delta^{1/2}$. Let $0<t<\delta/2$. By computation, one gets that for any $\beta \in \bR$, $2 \le p \le \infty$, and $t>0$,
    \[ \|\varphi\|_{L^2((t/2,t) \times B)} \lesssim_{p,\delta} t^\beta \| \I_{(t/2,t)} \varphi\|_{Z^p_\beta}. \]
    Here, $\I_{(t/2,t)} \varphi$ denotes the function $(s,y) \mapsto \I_{(t/2,t)}(s) \varphi(s,y)$. Applying this inequality to $u$ and $\nabla u$ yields that
    \[ \|u\|_{L^2((t/2,t);W^{1,2}(B))} \lesssim_{p,\delta} t^{1/2} \|\I_{(t/2,t)} u\|_{Z^p_{1/2}} + \|\I_{(t/2,t)} \nabla u\|_{Z^p_{0}}. \]
    Since $\partial_t u = \Div(A\nabla u + F)$, we get
    \[ \|\partial_t u\|_{L^2((t/2,t);W^{-1,2}(B))} \lesssim_{p,\delta} \|\I_{(t/2,t)} \nabla u\|_{Z^p_{0}} + \|\I_{(t/2,t)} F\|_{Z^p_{0}}. \]
    Hence, Lions' embedding theorem asserts $u \in C([t/2,t];L^2(B))$ with
    \[ \|u(t)\|_{L^2(B)} \lesssim_{p,\delta} t^{1/2} \|\I_{(t/2,t)} u\|_{Z^p_{1/2}} + \|\I_{(t/2,t)} \nabla u\|_{Z^p_{0}} + \|\I_{(t/2,t)} F\|_{Z^p_{0}}. \]
    Since the constant is independent of the center of $B$, we take supremum over all such balls $B$ to obtain
    \[ \|u(t)\|_{E^\infty_\delta} \lesssim_{p,\delta} t^{1/2} \|\I_{(t/2,t)} u\|_{Z^p_{1/2}} + \|\I_{(t/2,t)} \nabla u\|_{Z^p_{0}} + \|\I_{(t/2,t)} F\|_{Z^p_{0}}, \]
    which tends to 0 by Lebesgue's dominated convergence theorem. This proves the trace in $E^\infty_\delta$ and hence concludes Theorem \ref{thm:lions}.
    
    Theorems \ref{thm:ic} and \ref{thm:hc} can also be established analogously. In Theorem \ref{thm:hc}, to prove bounded extension of $\cE_A$ at the extreme point $\beta=0$ and $p=\infty$, one uses \eqref{e:heat-besov-eq-s=1-p=infty}.
    
    Next, we prove Theorem \ref{thm:unique}. For $p<\infty$, it readily follows from Sobolev embedding \eqref{e:sobolev-embed-T-Z}. For $p=\infty$, we first establish Proposition \ref{prop:trace} (with the canonical modification for $\beta=0$), using the same arguments in \cite[Propositions 6.2 and 9.1]{Auscher-Hou2024-HCL}. Then following the proof in Section \ref{ssec:pf-unique}, we may assume $u=\cR^{\bI}_{1/2}((A-\bI)\nabla u)$ (cf. \eqref{e:unique-u=RI1/2}) and then use the homotopy identity \eqref{e:homo-id} with the limits
    \[ \lim_{s \to 0} u(s) = 0 ~ \text{ in } E^\infty_\delta, \quad \lim_{s \to 0} \Gamma_A(t,s)^\ast h = \Gamma_A(t,0)^\ast h ~ \text{ in } E^1_\delta. \]
    This proves Theorem \ref{thm:unique} as desired.

    To prove Theorem \ref{thm:unique-u-Tpbeta+1}, we also use embedding of tent space for $p<\infty$. For $p=\infty$, we apply real interpolation to \eqref{e:unique-Tpbeta+1-u(s)Eps/16} to get
    \[ \|u(s)\|_{E^p_{s/16}} \lesssim s^{\beta+1/2} \|u\|_{Z^p_{\beta+1}}. \]
    Then the uniqueness follows from the same arguments in Section \ref{ssec:pf-unique-Tpbeta+1}.

    Theorem \ref{thm:rep} is a direct consequence of Theorem \ref{thm:unique} and Proposition \ref{prop:trace}. This completes the proof.
\end{proof}

It is worth noting that Theorem \ref{thm:besov} is based on our specific choice of realizations (in $\scrS'$) of homogeneous Besov spaces defined on $\scrS'/\scrP$. For instance, there is another family of realizations given by \cite[Definition 2.15]{Bahouri-Chemin-Danchin2011-book}, which does not contain polynomials for any $s$ and $p$, denoted by $\tilB^s_{p,p}$. The identification $\tilB^s_{p,p}+\bC = \DotB^s_{p,p}+\bC$ holds when $s<n/p$, but is not clear otherwise (see \cite[Remark 2.26]{Bahouri-Chemin-Danchin2011-book}).

\subsection*{Copyright}
A CC-BY 4.0 \url{https://creativecommons.
org/licenses/by/4.0/} public copyright license has been applied by the authors to the present document and will be applied to all subsequent versions up to the Author Accepted Manuscript arising from this submission.

\bibliographystyle{alpha}
\bibliography{sample}

\end{document}